\documentclass[11pt]{article}
\usepackage[top=1in, bottom=1in, left=1in, right=1in]{geometry}
\usepackage{tikz}
\usepackage{authblk}
\usepackage{amsmath,amsfonts,amssymb,amscd,amsthm,mathrsfs,float}
\usepackage{extarrows,textcomp}
\usepackage{graphicx,subfigure,epstopdf}
\usepackage{diagbox,multirow}
\usepackage{algorithm,algorithmic}
\usepackage{caption}
\usepackage{bm}
\usepackage{lineno}
\usepackage{xcolor}
\usepackage{makeidx,hyperref}
\usepackage{bm}
\usepackage{cases}

\theoremstyle{definition}
\newtheorem{remark}{Remark}{}
\numberwithin{equation}{section}

\DeclareMathOperator*{\argmax}{arg\,max}
\DeclareMathOperator*{\argmin}{arg\,min}

\newcommand{\Loss}{\mathrm{Loss}}

\newcommand{\transpose}{\mathsf{T}}

\makeatletter
\newcommand*{\extendadd}{
  \mathbin{
    \mathpalette\extend@add{}
  }
}
\newcommand*{\extend@add}[2]{
  \ooalign{
    $\m@th#1\leftrightarrow$%
    \vphantom{$\m@th#1\updownarrow$}
    \cr
    \hfil$\m@th#1\updownarrow$\hfil
  }
}
\makeatother

\title{Adaptive neural network basis methods for partial differential equations with low-regular solutions}
\date{}

\author[1]{Jianguo~Huang
\thanks{E-mail: jghuang@sjtu.edu.cn. The work of this author was partially supported by the National
Natural Science Foundation of China (Grant No. 12071289), and the Strategic Priority Research Program of the Chinese Academy of Sciences (Grant No. XDA25010402).}
}
\author[1]{Haohao~Wu
\thanks{E-mail: wu1150132305@sjtu.edu.cn. Corresponding author.}
}
\author[2]{Tao~Zhou
\thanks{tzhou@lsec.cc.ac.cn. The work of this author was partially supported by the National
Natural Science Foundation of China (Grant No. 12288201) and the Youth Innovation Promotion Association (CAS). }}

\affil[1]{School of Mathematical Sciences, and MOE-LSC, Shanghai Jiao Tong University, Shanghai 200240, China.
}
\affil[2]{Institute of Computational Mathematics and Scientific/Engineering Computing, Academy of Mathematics and Systems Science, Chinese Academy of Sciences, Beijing, China.}

\begin{document}

\maketitle
\begin{abstract}
This paper aims to devise an adaptive neural network basis method for numerically solving a second-order semilinear partial differential equation (PDE) with low-regular solutions in two/three dimensions. The method is obtained by combining basis functions from a class of shallow neural networks and the resulting multi-scale analogues, a residual strategy in adaptive methods and the non-overlapping domain decomposition method. At the beginning, in view of the solution residual, we partition the total domain $\Omega$ into $K+1$ non-overlapping subdomains, denoted respectively as $\{\Omega_k\}_{k=0}^K$, where the exact solution is smooth on subdomain $\Omega_{0}$ and low-regular on subdomain $\Omega_{k}$ ($1\le k\le K$). Secondly, the low-regular solutions on different subdomains \(\Omega_{k}\)~($1\le k\le K$) are approximated by neural networks with different scales, while the smooth solution on subdomain \(\Omega_0\) is approximated by the initialized neural network.
Thirdly, we determine the undetermined coefficients by solving the linear least squares problems directly or the nonlinear least squares problem via the Gauss-Newton method. The proposed method can be extended to multi-level case naturally. Finally, we use this adaptive method for several peak problems in two/three dimensions to show its high-efficient computational performance.
\end{abstract}

\section{Introduction}
In the past few years, the machine learning including neural networks and artificial intelligence (AI) has experienced rapid and significant advances in computer science and data analysis (cf. \cite{2015-DL,2016-DL}). This technology has represented a cornerstone of innovation across various industries, continues to transform how we live and work from enhancing everyday experiences to driving groundbreaking discoveries. More recently, the machine learning method has also become an important approach to numerically solving partial differential equations (PDEs), which is a heart topic in computational and applied mathematics. Initial exploration of such study can be traced back to significant works in the 1990s (cf.   \cite{1994-Dissanayake,1998-Lagaris}). The typical methods along this line include but not limited to the following methods. The Deep Ritz method was introduced in \cite{2018-DRM} for solving elliptic problems by combining the Ritz method and deep neural works (DNNs). The weak adversarial  networks \cite{2020-WAN} is a machine learning method for solving a weak solution of a PDE based on its weak formulation and using the solution ansatz by DNNs (see also \cite{2023-Friedrichs-learning-Huang}). Another important classes of methods is the so-called Physics-Informed Neural Network (PINN) methods and their variants; we refer the reader to \cite{2019-PINN,2021-DeepXDE-Karniadakis,2019-fPINNs-Karniadakis,2020-PINN-Analysis-Karniadakis,2020-cPINNs-Karniadakis,2020-XPINNs-Karniadakis,2021-hp-Vpinns-Karniadakis} for details. The applications of such methods in computational mathematics can be found in (cf.\cite{2020-E}). In all of the above machine learning methods, DNNs are used to parameterize the solution of PDEs and the related parameters are identified by minimizing an optimization problem formulated from the PDEs. The remarkable advantage of such methods connecting with DNNs is that they can overcome the so-called ``curse of dimensionality" (cf. \cite{E-2018-HD,2019-HD-E,2020-HD-BSPDE,2020-HD-heat-equation,2021-HD-KPDEs}). Furthermore, to improve computational efficiency, adaptive sampling techniques are combined with the PINN approach for solving PDEs with low-regular solutions (cf.\cite{2021-DeepXDE-Karniadakis,2021-DAS,2023-residual-based-as-Lulu,2023-FIPINN-zhou1,2023-FIPINN-zhou2}). Moreover, such techniques are also combined with the deep Ritz method for solving PDEs in \cite{2023-adaptive-sample-DRM}.



Despite the DNN-based machine learning method can get over "curse of dimensionality", when solving PDEs in low dimensions, one is tempted to use shallow neural network as the solution asatz to balance computational accuracy and cost. Moreover, the weight/bias coefficients in the hidden layer are pre-set to random values and fixed, while the training parameters consist of the weight coefficients of output layer and they are obtained by using existing linear or nonlinear least squares solvers. In other words, one can use the shallow network to produce a basis of the underlying admissible space of the infinite-dimensional minimization problem associated with a PDE under discussion. An important advantage of such a basis is that it is mesh-free. We mention that randomness has been used in neural networks for a long time (cf. \cite{2017-RNN}). Randomized neural networks can be traced back to Turing's unorganized machine and Rosenblatt's perception \cite{2011-Turing,1958-perceptron}.  Later on, extreme Learning machine (ELM) was proposed in \cite{2006-ELM,2011-ELM,2015-ELM} for solving linear classification or regression problems by combining the previous random neural network with the linear least squares method. This method was also used for solving  ordinary differential equations (cf. \cite{2018-ELM-Legendre-ODE,2020-ELM-Legendre-PDE}) and partial differential equations (PDEs) (cf. \cite{2019-ELM-PDE,2020-PIELM}) . The ELM method was further combined with the non-overlapping domain decomposition technique to solve  PDEs effectively (cf. \cite{2021-ELM-Dong-Li,2021-LoELM-Dong-Li}). In \cite{2024-ELM-HD-Dong-Wang}, the method was even applied to solve high-dimensional partial differential equations. The combination of the random neural network and the partition of unity technique gives rise to the random feature method (cf. \cite{2022-chen-RFM}). More recently, transferable neural networks (cf. \cite{2024-Ju}) generate hidden layer parameters by re-parameterizing the hidden neurons and using auxiliary functions to tune. It is proved that the hidden layer neurons are uniformly distributed within the unit cell. Numerical results indicate that the computational accuracy is high if the solution is smooth enough. However, the solution accuracy would deteriorate remarkably if the exact solution has low regularity.

In this work, we are going to make use of the technique in constructing the basis functions from a shallow neural network, the non-overlapping domain decomposition (DDM), and the multi-scale neural networks \cite{2020-Mscale-Xu} to produce an adaptive neural network basis (ANNB) method for solving second-order PDEs with low-regular solutions. The main idea of the ANNB method can be summarized in three steps. First of all, in view of the solution residual,  we partition the total domain $\Omega$ into $K+1$ non-overlapping subdomains, denoted respectively as $\Omega_{0},\Omega_{1},\cdots,\Omega_{K}$, where the solution of the target PDE is smooth on subdomain $\Omega_{0}$ and is low-regular on each subdomain $\Omega_{k}$ ($1\le k\le K$).
Secondly, the low-regular solutions on different subdomains \(\Omega_{k}\)~($1\le k\le K$) are approximated by neural networks with different scales, while the smooth solution on subdomain \(\Omega_0\) is approximated by the initialized neural network. Thirdly, we determine the unknown coefficients  by solving the linear least squares problems or the nonlinear least squares problems in view of the Gauss-Newton method. It is worth noting that  the domain decomposition can be multi-leveled, though we only consider the two-level case for simplicity. Finally, we use this adaptive method to solve peak problems in two and three dimensions, which show the resulting numerical solution is high accurate.

The rest of the paper is organized as follows. In section 2, we show how to construct neural network basis functions and a DDM-based neural network basis method (NNBM) for second-order linear and nonlinear PDEs. In section 3, we present the ANNB  method in detail  for solving  second-order PDEs with low-regular solutions. In section 4, several numerical tests are provided to demonstrate the effectiveness of the ANNB method. Some conclusions are given in section 5. Finally, an appendix is provided to show the implementation details of our ANNM method.

\section{Construction of neural network basis functions and a DDM-based NNBM  for second order linear and nonlinear PDEs}\label{sec:TNNs}
Let $\Omega \subset \mathbb{R}^{d}$ ($d=1,2,3$) be a bounded domain with the Lipschitz boundary $\partial\Omega$. Consider the following generic second-order linear or nonlinear boundary value problem
\begin{equation}\label{pro:ini}
    \begin{cases}
    \mathcal{L}u(\bm{x}) = f(\bm{x}) & \text{for } \bm{x} \in \Omega, \\
    u(\bm{x}) = g(\bm{x})  & \text{for } \bm{x} \in \partial \Omega,
    \end{cases}
\end{equation}
where $\bm{x}:=\left(x_1, \ldots, x_{d}\right)^{\intercal}$ are independent variables,  $\mathcal{L}$ is a second-order semi-linear differential operator from Sobolev space $H^2(\Omega)$ to the  space $L^{2}(\Omega)$, $g$ is the restriction to $\partial\Omega$ of a certain function in $H^2(\Omega)$ and the scalar function $u$ denotes the solution of the above problem.
\subsection{Construction of neural network basis}
As is well-known, we can express a shallow neural network in the form
\begin{equation}\label{eq:NN}
 u_{\mathrm{NN}}(\bm{x}) = \sum_{m=1}^M \alpha_m \sigma\left(\bm{w}^{\transpose}_m \bm{x} + b_m\right) + \alpha_0,
\end{equation}
where $M$ represents the number of hidden neurons, the row vector $\bm{w}_m = \left(w_{m, 1}, \ldots, w_{m, d}\right)^{\transpose}$ means the weights of the $m$-th hidden neuron, the scalar $b_m$ means the bias of the $m$-th hidden neuron, the row vector $\bm{\alpha} = \left(\alpha_0, \alpha_1, \ldots, \alpha_M\right)^{\transpose}$ includes the weights and bias of the output layer, and $\sigma(\cdot)$ denotes the activation function. Besides, we assume the weights/bias coefficients in the hidden layer are pre-set during the computation. In other words, $u_{\mathrm{NN}}(\bm{x})$ lies in the admissible space
\[
V_M={\rm span}\{\psi_0,\psi_1,\cdots, \psi_M\},
\]
where
\[
\psi_{0}(\bm{x})=1, \psi_{m}(\bm{x})=\sigma({\bm{\omega}_{m}}^{\transpose}\bm{x}+b_{m})~(1\le m \le M),
\]
are $M+1$ neural network basis functions.

There are two strategies to pre-set the weights/bias coefficients in the hidden layer.

The first strategy comes from \cite{2021-LoELM-Dong-Li,2022-chen-RFM,2022-ELM-Dong-Yang}. Without loss of generality, assume $\Omega=[-1,1]^d$, the weight/bias coefficients $\{\bm{w}_{m}\}_{m=1}^{M}$ and $\{b_{m}\}_{m=1}^{M}$ are generated  from a uniform distribution on $[-R,R]^d$ and on $[-R,R]$, where $R$ is a user-defined constant parameter.

The second strategy is given by \cite{2024-Ju}.
In the unit ball domain, i.e., $B_{1}(\bm{0})=\{\bm{x},\|\bm{x}\|_{2}\le 1\}\subset\mathbb{R}^{d}$, the weights/bias coefficients  are generated by the following three steps for the activation function  $\tanh (\cdot)$.

\begin{enumerate}
\item Firstly, each hidden
neuron can be rewritten as
\begin{equation}\label{eq:24-Ju-1}
\sigma(\bm{w}_{m}^{\transpose}\bm{x}+b_{m})=\sigma(\gamma_{m}(\bm{a}_{m}^{\transpose}\bm{x}+r_{m}))
\end{equation}
where $\|\bm{a}_{m}\|_{2}=1$ and $\gamma_{m}\ge 0$. Then the parameter $(\bm{\omega}_{m},b_{m})$ is re-parameterized into the
location parameter $(\bm{a}_{m},r_{m})$ and shape parameter $\gamma_{m}$. Their equivalence is given by
\begin{equation*}
\begin{cases}
 \bm{w}_{m} =\gamma_m \bm{a}_{m},\\
b_m = \gamma_m r_m.
\end{cases}
\Longleftrightarrow
\begin{cases}
 \bm{a}_{m} =\frac{\bm{w}_{m}}{\|\bm{w}_{m}\|_{2}},\\
r_m =\frac{b_{m}}{\|\bm{w}_{m}\|_{2}},\\
\gamma_{m}=\|\bm{w}_{m}\|_{2}.
\end{cases}
\end{equation*}
    \item The second step is to determine the location parameters $\left\{\left(\bm{a}_m, r_m\right)\right\}_{m=1}^M$ in Eq.\eqref{eq:24-Ju-1} as follows:
    \begin{equation}\label{eq:24-Ju-2}
    \bm{a}_m=\frac{\bm{X}_m}{\left\|\bm{X}_m\right\|_2}, r_m=U_m, m=1, \ldots, M,
    \end{equation}
where $\bm{X}_{m}$  are i.i.d $d$-dimensional standard Gaussian distribution, $U_{m}$ follow i.i.d. uniform distribution
on $[0, 1]$. This approach
yields a set of uniformly
distributed partition hyperplanes in the unit ball $B_{1}(\bm{0})$.
\item  The third step is tune the shape parameters  $\left\{\gamma_m\right\}_{m=1}^M$ in Eq.\eqref{eq:24-Ju-1}. For simplicity, let all neurons share the same shape parameter value, i.e., $\gamma=\gamma_m$ for $m=1, \ldots, M$.
   The main idea is to use realizations of Gaussian random fields as the auxiliary functions for tuning  the shape parameter $\gamma$ without using any information about a PDE. Please refer to  Algorithm 1 in \cite{2024-Ju} for details.
\end{enumerate}

\subsection{A DDM-based NNBM}\label{sec:DDBTLM}
We partition the domain $\Omega$ into $K+1~(K\ge 0)$ non-overlapping sub-domains,
$$
\bar{\Omega}=\cup_{k=0}^{K}\bar{\Omega}_{k},~\bar{\Omega}_{k_{1}}\cap\bar{\Omega}_{k_{2}}=\emptyset,~1\le k_{1}\ne k_{2}\le K,
$$
where $\Omega_{k}$ denotes the $k$-th subdomain. More specifically, the sub-domain
 $\Omega_{k}~(1\le k\le K)$ only share a common boundary with sub-domain $\Omega_{0}$ and we will denote this common boundary by $\Gamma_{k}$, i.e., $\Gamma_{k}= \partial\Omega_{k}\cap\partial{\Omega}_{0}$.
 Additionally, the boundary \(\partial \Omega_{k}\) also includes the portion that intersects with the boundary of the original domain \(\partial \Omega\).
It is obvious that  $\Omega_{0}=\Omega$ when $K=0$.

We enforce the $C^{1}$ continuity conditions on the interface $\Gamma_{k}~(1\le k\le K)$, and then  the  problem \eqref{pro:ini} can be rewritten equivalently as the following $(K+1)$-subdomain problem :
\begin{equation}\label{pro:MD}
    \begin{cases}
    \mathcal{L}u_{0}(\bm{x})=f(\bm{x}) &\text{for} \bm{x} \in \Omega_{0}, \\
    u_{0}(\bm{x})=g(\bm{x})  &\text{for} \bm{x} \in \partial\Omega_{0} \cap\partial\Omega, \\
    \mathcal{L}u_{1}(\bm{x})=f(\bm{x}) &\text{for} \bm{x} \in \Omega_{1}, \\
    u_{1}(\bm{x})=g(\bm{x})  &\text{for} \bm{x} \in \partial\Omega_{1} \cap\partial\Omega, \\
    u_{1}(\bm{x})=u_{0}(\bm{x}) &\text{for} \bm{x} \in  \Gamma_{1},\\
    \frac{\partial u_{1}(\bm{x})}{\partial \bm{n}_{1}}=\frac{\partial u_{0}(\bm{x})}{\partial \bm{n}_{1}} &\text{for} \bm{x} \in  \Gamma_{1},
    \\
    \qquad\vdots&
    \\
    \mathcal{L}u_{K}(\bm{x})=f(\bm{x}) &\text{for} \bm{x} \in \Omega_{K}, \\

    u_{K}(\bm{x})=g(\bm{x})  &\text{for} \bm{x} \in \partial\Omega_{K} \cap\partial\Omega,\\

    u_{K}(\bm{x})=u_{0}(\bm{x}) &\text{for} \bm{x} \in  \Gamma_{K},\\
    \frac{\partial u_{K}(\bm{x})}{\partial \bm{n}_{K}}=\frac{\partial u_{0}(\bm{x})}{\partial \bm{n}_{K}} &\text{for} \bm{x} \in  \Gamma_{K},
    \\
    \end{cases}
\end{equation}
where $\bm{n}_{k} $ is the outward normal vector along the interface $\Gamma_{k}$. It is worth noting that if $K=0$, problem \eqref{pro:MD} is exactly problem \eqref{pro:ini}.

In each sub-domain $\Omega_{k}~(0\le k\le K)$, we represent the solution of $u(\bm{x})$ by using a set of neural network basis functions which span the  admissible space
$V_{M,k}=\mathrm{span} \{\psi_{0,k},\psi_{1,k},\cdots,\psi_{M,k}\}$, i.e.,
\begin{equation*}
\bm{u}(x)\approx \Tilde{u}_{k}(\bm{\alpha}_{k},\bm{x}) = \sum_{m=0}^{M_{k}} \alpha_{m,k} \psi_{m,k}(\bm{x}) ~ \text{for}~ \bm{x}\in \Omega_{k}.
\end{equation*}
Here, $ M_{k} $ is the number of neural network basis functions, $ \{\psi_{m,k}(\bm{x})\}_{m=0}^{M_{k}}$ denote the neural network basis functions, and $\bm{\alpha}_{k} = (\alpha_{0,k}, \cdots, \alpha_{M_{k},k})^{\transpose}$ are the trainable parameters on $ \Omega_{k}$.

Corresponding to  \eqref{pro:MD}, we have $3K+2$ sets of collocation points: $\bm{X}_{f_{k}}~(0\le k\le K)$, the set of interior collocation points in $\Omega_{k}$, $\bm{X}_{g_{k}}~(0\le k\le K)$, the set of boundary collocation points on $\Omega_{k}\cap\partial\Omega$ and $\bm{X}_{\Gamma_k}~(1\le k\le K)$, the set of interface collocation points on $\Gamma_{k}$.

The loss function is built on the strong formulation of \eqref{pro:MD} at collocation points so as to produce the following least square problem:

\begin{equation}\label{pro:MD-min}
  \begin{aligned}
    \min_{\bm{\alpha}^{0}, \cdots, \bm{\alpha}^{K}}&\left\{ \sum_{k=0}^{K}\left(\sum_{\bm{x}_{f_{k}}\in \bm{X}_{f_{k}}}|\mathcal{L}\tilde{u}_{k}(\bm{x}_{f_{k}})-f(\bm{x}_{f_{k}})|^2+\sum_{\bm{x}_{g_{k}}\in \bm{X}_{g_{k}}}|\tilde{u}_{k}(\bm{x}_{g_{k}})-g(\bm{x}_{g_{k}})|^2\right)\right.\\
    &\left.+\sum_{k=1}^{K}\sum_{\bm{x}_{\Gamma_{k}}\in \bm{X}_{\Gamma_{k}}}\left(\left|\tilde{u}_{k}(\bm{x}_{\Gamma_{k}})-\tilde{u}_{0}(\bm{x}_{\Gamma_{k}})\right|^2+\left|\frac{\partial \tilde{u}_{k}(\bm{x}_{\Gamma_{k}})}{\partial \bm{n}_{k}}-\frac{\partial \tilde{u}_{0}(\bm{x}_{\Gamma_k})}{\partial \bm{n}_{k}}\right|^2\right)\right\}.
  \end{aligned}
\end{equation}

When the differential operator  $\mathcal{L}$ is linear, the above problem can be reformulated as a linear least square problem in the form
\begin{equation}\label{pro:MD-min-linear}
\begin{aligned}
\min_{\bm{\alpha}^{0}, \cdots, \bm{\alpha}^{K}} &\left\{
\sum_{k=0}^{K}\left(\sum_{\bm{x}_{f_{k}}\in \bm{X}_{f_{k}}}\left|\sum_{m=0}^{M_{k}} \alpha_{m,k}\mathcal{L}\psi_{m,k}(\bm{x}_{f_{k}})-f(\bm{x}_{f_{k}})\right|^2\right.\right.\\
&\left. +\sum_{\bm{x}_{g_{k}}\in \bm{X}_{g_{k}}}\left|\sum_{m=0}^{M_{k}} \alpha_{m,k}\psi_{m,k}(\bm{x}_{g_{k}})-g(\bm{x}_{g_{k}})\right|^2\right)\\
&+\sum_{k=1}^{K}\sum_{\bm{x}_{\Gamma_{k}}\in \bm{X}_{\Gamma_{k}}}\left(\left| \sum_{m=0}^{M_{k}} \alpha_{m,k} \psi_{m,k}(\bm{x}_{\Gamma_{k}})- \sum_{m=0}^{M_{0}} \alpha_{m,0} \psi_{m,k}(\bm{x}_{\Gamma_{k}})\right|^2\right.\\
&\left.\left.+\left|\frac{\sum_{m=0}^{M_{k}}\alpha_{m,k}\partial \psi_{m,k}(\bm{x}_{\Gamma_{k}})}{\partial \bm{n}_{k}}-\frac{\sum_{m=0}^{M_{0}}\partial\alpha_{m,0}\ \psi_{m,0}(\bm{x}_{\Gamma_{k}})}{\partial \bm{n}_{k}}\right|^2\right)\right\},
\end{aligned}
\end{equation}
which can be solved directly.


When the differential operator $\mathcal{L}$ is nonlinear, we employ the Gauss-Newton method to solve the least squares problem \eqref{pro:MD-min} (see also \cite{2024-Huang-Wu-GN}). Concretely speaking, at the $n$-th iteration step, we denote $u_{k}^{n}(\bm{x})~(0\le k\le K)$ as the approximation of the solution $u_{k}(\bm{x})$, and $v^{n}_{k}(\bm{x})$ as the increment field that needs to be computed. Then
$$
u_{k}^{n+1}(\bm{x})= u_{k}^{n}(\bm{x})+v_{k}^{n}(\bm{x}).
$$

The increment $v_{k}^{n}(\bm{x})$ is determined by the linearized version of \eqref{pro:MD} given below
\begin{equation}\label{pro:GN}
        \begin{cases}
    D\mathcal{L}(u_{0}^{n};v_{0}^{n})(\bm{x})  = f(\bm{x})-\mathcal{L}u_{0}^{n}(\bm{x})  &\text{for} \bm{x} \in \Omega_{0}, \\

    v_{0}^{n}(\bm{x}) = g(\bm{x}) -u_{0}^{n}(\bm{x}) &\text{for} \bm{x} \in \partial\Omega_{0} \cap\partial\Omega,\\
    D\mathcal{L}(u_{1}^{n};v_{1}^{n}) (\bm{x}) = f(\bm{x})-\mathcal{L}u_{1}^{n}(\bm{x}) &\text{for} \bm{x} \in \Omega_{1}, \\
    v_{1}^{n}(\bm{x}) = g(\bm{x}) -u_{1}^{n}(\bm{x}) &\text{for} \bm{x} \in \partial\Omega_{1} \cap\partial\Omega,\\
    v_{1}^{n}(\bm{x})-v_{0}^{n}(\bm{x})=u_{0}^{n}(\bm{x})-u_{1}^{n}(\bm{x}) &\text{for} \bm{x} \in  \Gamma_{1},\\
    \frac{\partial v_{1}^{n}(\bm{x})}{\partial \bm{n}_{1}}-\frac{\partial v_{0}^{n}(\bm{x})}{\partial \bm{n}_{1}}=\frac{\partial u_{1}^{n}(\bm{x})}{\partial \bm{n}_{1}}-\frac{\partial u_{0}^{n}(\bm{x})}{\partial \bm{n}_{1}} &\text{for} \bm{x} \in  \Gamma_{1},
    \\
    \qquad\qquad\vdots\\
     D\mathcal{L}(u_{K}^{n};v_{K}^{n})(\bm{x})  = f(\bm{x})-\mathcal{L}u_{k}^{n}(\bm{x}) &\text{for} \bm{x} \in \Omega_{K}, \\
    v_{K}^{n}(\bm{x}) = g(\bm{x}) -u_{k}^{n}(\bm{x}) &\text{for} \bm{x} \in \partial\Omega_{K} \cap\partial\Omega,\\

    v_{K}^{n}(\bm{x})-v_{0}^{n}(\bm{x})=u_{0}^{n}(\bm{x})-u_{K}^{n}(\bm{x}) &\text{for} \bm{x} \in  \Gamma_{K},\\
    \frac{\partial v_{K}^{n}(\bm{x})}{\partial \bm{n}_{K}}-\frac{\partial v_{0}^{n}(\bm{x})}{\partial \bm{n}_{K}}=\frac{\partial u_{K}^{n}(\bm{x})}{\partial \bm{n}_{K}}-\frac{\partial u_{0}^{n}(\bm{x})}{\partial \bm{n}_{K}} &\text{for} \bm{x} \in  \Gamma_{K},
    \\
    \end{cases}
\end{equation}
where $D\mathcal{L}(u(\bm{x});v(\bm{x}))$ is the G\^{a}teaux derivative (cf. \cite{Nonlinear-Analysis-05}) of $\mathcal{L}$ at $u(\bm{x})$ in the direction of $v(\bm{x})$.

Then we represent $u_{k}^{n}(\bm{x})$ and $v_{k}^{n}(\bm{x})$ by using a set of neural network basis functions which span the  admissible space
$V_{M,k}=\mathrm{span} \{\psi_{0,k},\psi_{1,k},\cdots,\psi_{M,k}\}$, i.e.,
$$
\bm{u}_{k}^{n}(x)\approx \Tilde{u}_{k}(\bm{\alpha_{k}^{n}},\bm{x}) = \sum_{m=0}^{M_{k}} \alpha_{m,k}^{n} \psi_{m,k}(\bm{x}),~
\bm{v}_{k}^{n}(x)\approx \Tilde{u}_{k}(\bm{a}_{k}^{n},\bm{x}) = \sum_{m=0}^{M_{k}} a_{m,k}^{n} \psi_{m,k}(\bm{x}) ~ \text{for}~ \bm{x}
\in\Omega_{k}.
$$
Define $\bm{\alpha}_{k}^{n}=(\alpha_{0,k}^{n},\alpha_{1,k}^{n},\cdots,\alpha_{M_{k},k}^{n})^{\transpose}$ and $\bm{a}_{k}^{n}=(a^{n}_{0,k},a^{n}_{1,k},\cdots,a^{n}_{M_{k},k})^\transpose$.

Following a similar derivation to obtain problem \eqref{pro:MD-min-linear}, we can transform the problem \eqref{pro:GN} as a linear least square problem given below.

\begin{equation}\label{pro:MD-min-nonlinear}
\begin{aligned}
\min_{\bm{a}_{0}^{n}, \cdots, \bm{a}_{K}^{n}} &\left\{
\sum_{k=0}^{K}\left(\sum_{\bm{x}_{f_{k}}\in \bm{X}_{f_{k}}}\left|\sum_{m=0}^{M_{k}} a_{m,k}^{n}D\mathcal{L}(\tilde{u}_{0}^{n};\psi_{m,k})(\bm{x}_{f_{k}})-f(\bm{x}_{f_{k}})+\mathcal{L}\tilde{u}_{0}^{n}(\bm{x}_{f_{k}})\right|^2 \right.\right.\\
\qquad \qquad &\left.+\sum_{\bm{x}_{g_{k}}\in \bm{X}_{g_{k}}}\left|\sum_{m=0}^{M_{k}} a_{m,k}^{n}\psi_{m,k}(\bm{x}_{g_{k}})-g(\bm{x}_{g_{k}})+\tilde{u}_{k}^{n}(\bm{x}_{g_{k}})\right|^2\right)+\\
&\sum_{k=1}^{K}\sum_{\bm{x}_{\Gamma_{k}}\in \bm{X}_{\Gamma_{k}}}\left(\left| \sum_{m=0}^{M_{k}} a_{m,k} \psi_{m,k}(\bm{x}_{\Gamma_{k}})- \sum_{m=0}^{M_{0}} a_{m,0} \psi_{m,k}(\bm{x}_{\Gamma_{k}})+\tilde{u}_{k}(\bm{\alpha}_{k}^{n},\bm{x}_{\Gamma_{k}})-\tilde{u}_{0}(\bm{\alpha}_{0}^{n},\bm{x}_{\Gamma_{k}})\right|^2\right.+\\
&\left.\left.\left|\frac{\sum_{m=0}^{M_{k}}\alpha_{m,k}\partial \psi_{m,k}(\bm{x}_{\Gamma_{k}})}{\partial \bm{n}_{k}}-\frac{\sum_{m=0}^{M_{0}}\alpha_{m,0}\ \psi_{m,0}(\bm{x}_{\Gamma_{k}})}{\partial \bm{n}_{k}}+\frac{\partial\tilde{u}_{k}(\bm{\alpha}_{k}^{n},\bm{x}_{\Gamma_{k}})}{\partial\bm{n}_{k}}-\frac{\partial\tilde{u}_{0}(\bm{\alpha}_{0}^{n},\bm{x}_{\Gamma_{k}})}{\partial\bm{n}_{k}}\right|^2\right)\right\}
\end{aligned}
\end{equation}
which can be solved directly.

The details  of  DDM-based NNBM for the second-order nonlinear problem \eqref{pro:MD-min} are outlined in Algorithm \ref{alg:DDM-nonlinear}.
\begin{algorithm}[H]
\caption{A DDM-based NNBM for the second-order  problem \eqref{pro:MD}}
\label{alg:DDM-nonlinear}
\begin{algorithmic}[1]
\STATE{$\mathbf{Input}$: $\left\{\psi_{m,k}(\bm{x})\right\}_{m=0}^{M_{k}}~(0\le k\le K)$,~neural network basis functions for sub-domain $\Omega_{k}$; $\bm{X}_{f_{k}}~(0\le k\le K)$, the set of interior collocation points in $\Omega_{k}$; $\bm{X}_{g_{k}}~(0\le k\le K)$ the set of boundary collocation points on $\partial\Omega_{k}\cap\partial\Omega_{k}$, ;  $\bm{X}_{\Gamma_{k}}~(1\le k\le K)$, the set of collocation points on $\Gamma_{k}$; $N$,the maximum number of iteration steps ; $tol$, tolerance.}
\STATE{$\mathbf{Output}:\bm{\alpha}_{0}^{*},\cdots,\bm{\alpha}_{K}^{*}$.}

\IF{$\mathcal{L}$ is linear operator}
\STATE{
Obtain $\bm{\alpha}^{*}=({\bm{\alpha}_{0}^{*}}^\transpose,\cdots,{\bm{a}_{K}^{*}}^\transpose)^\transpose$ by solving the following linear least square problem:
  $$
  \min_{\bm{\alpha}}\|\bm{F}\bm{\alpha}-\bm{T}\|_{2}^{2}.
  $$
  where  the matrices \(\bm{F}\) and \(\bm{T}\)  are constructed  according to the problem \eqref{pro:MD-min-linear} and $\bm{\alpha}=({\bm{a}_{0}}^\transpose,{\bm{a}_{1}}^\transpose,\cdots,{\bm{a}^{K}}^\transpose)^\transpose$.}

\ELSE

 \STATE{Initial parameters $\bm{\alpha}_{k}^{0}=(\alpha^0_{0,k},\alpha^0_{1,k},\cdots,\alpha^0_{M_{k},k})^{\transpose}=(0,0,\cdots,0)^{\transpose}~(0\le k\le K),~n=0$.}
 \WHILE{$n\le N$}

\STATE{
Obtain $\bm{a}^{n,*}=({\bm{a}_{0}^{n,*}}^\transpose,\cdots,{\bm{a}_{K}^{n,*}}^\transpose)^\transpose$ by solving the following linear least square problem:
$$
\min_{\bm{a}^{n}}\|\bm{F}^{n}\bm{a}^{n}-\bm{T}^{n}\|_{2}^{2}.
$$
where  the matrices \(\bm{F}^{n}\) and \(\bm{T}^{n}\)  are constructed  according to the problem \eqref{pro:MD-min-nonlinear} and $\bm{a}^{n}=({\bm{a}_{0}^{n}}^\transpose,{\bm{a}_{1}^{n}}^\transpose,\cdots,{\bm{a}^{K}_{n}}^\transpose)^\transpose$.}
 \IF{$n\ge 1$}
\STATE{$Re_{mse}=\frac{\left|\|\bm{F}^{n}\bm{a}^{n,*}-\bm{T}^{n}\|_{2}^{2}-\|\bm{F}^{n-1}\bm{a}^{n-1,*}-\bm{T}^{n-1}\|_{2}^{2}\right|}{\|\bm{F}^{n-1}\bm{a}^{n-1,*}-\bm{T}^{n-1}\|_{2}^{2}}$}
\IF{$Re_{mse}<tol$}
    \STATE{Break.}
    \ENDIF
\ENDIF
\STATE{
$\bm{\alpha}^{n+1}_{k}=\bm{\alpha}_{k}^{n}+\bm{a}_{K}^{n,*}~(0\le k\le K)$.}
 \ENDWHILE
\STATE{$\bm{\alpha}^{*}_{k}=\bm{\alpha}_{k}^{n+1}~(0\le k\le K)$.}
\ENDIF
\end{algorithmic}
\end{algorithm}

\section{The ANNB method}

\begin{figure}[H]
  \begin{minipage}{0.32\textwidth}
    \centering
    \includegraphics[width=\textwidth,height=5cm]{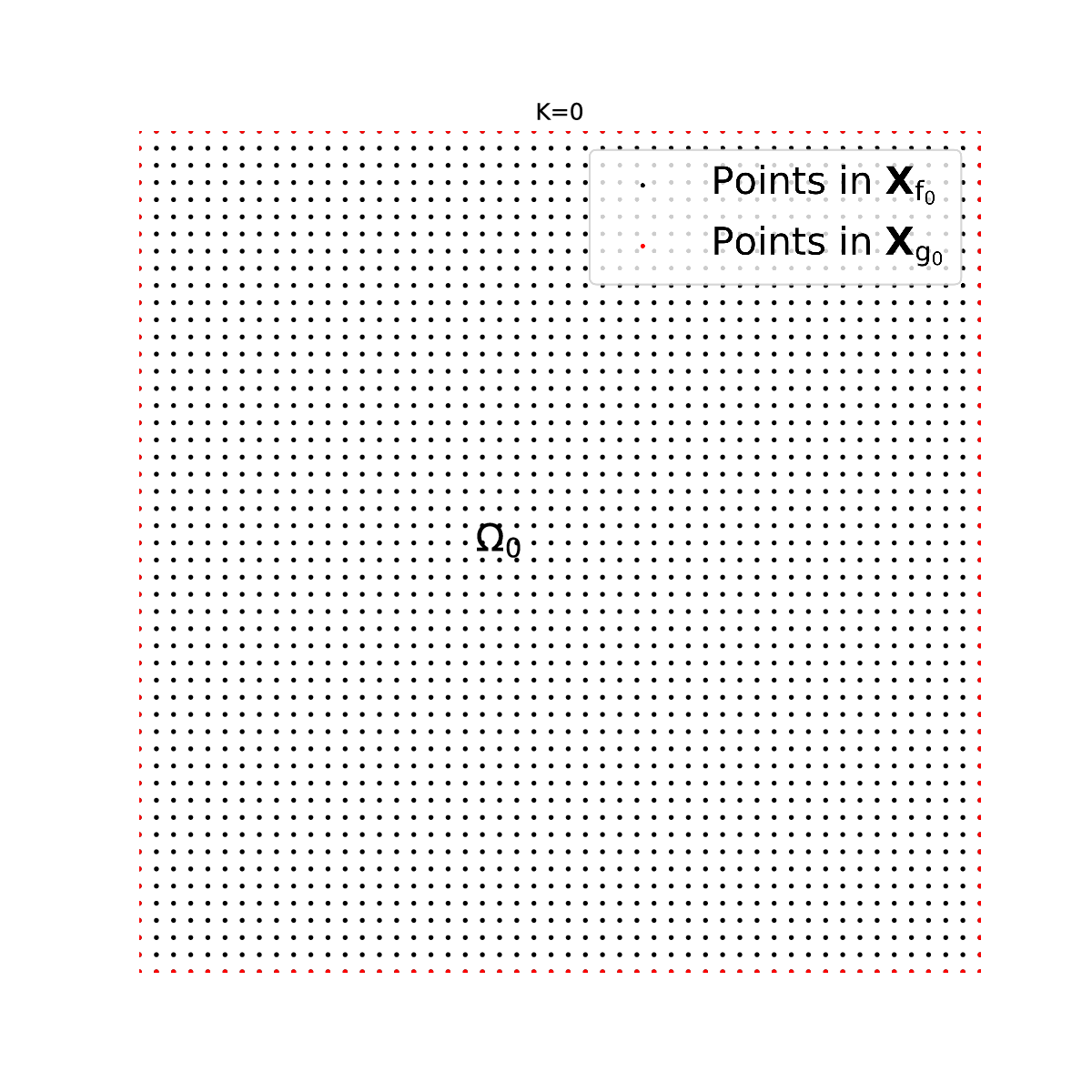}
    \end{minipage}
  \begin{minipage}{0.32\textwidth}
    \centering
    \includegraphics[width=\textwidth,height=5cm]{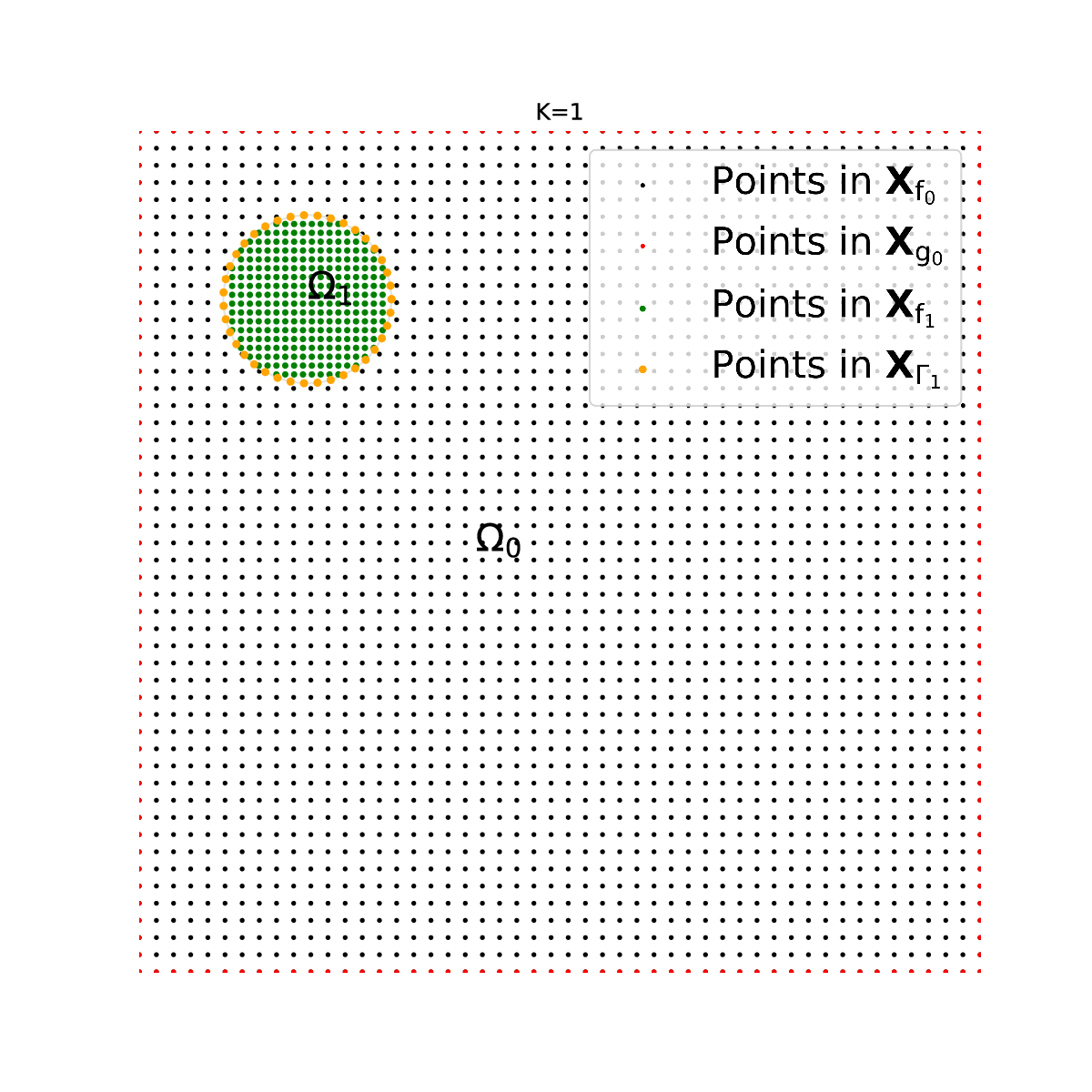}
  \end{minipage}
\begin{minipage}{0.32\textwidth}
  \includegraphics[width=\textwidth,height=5cm]{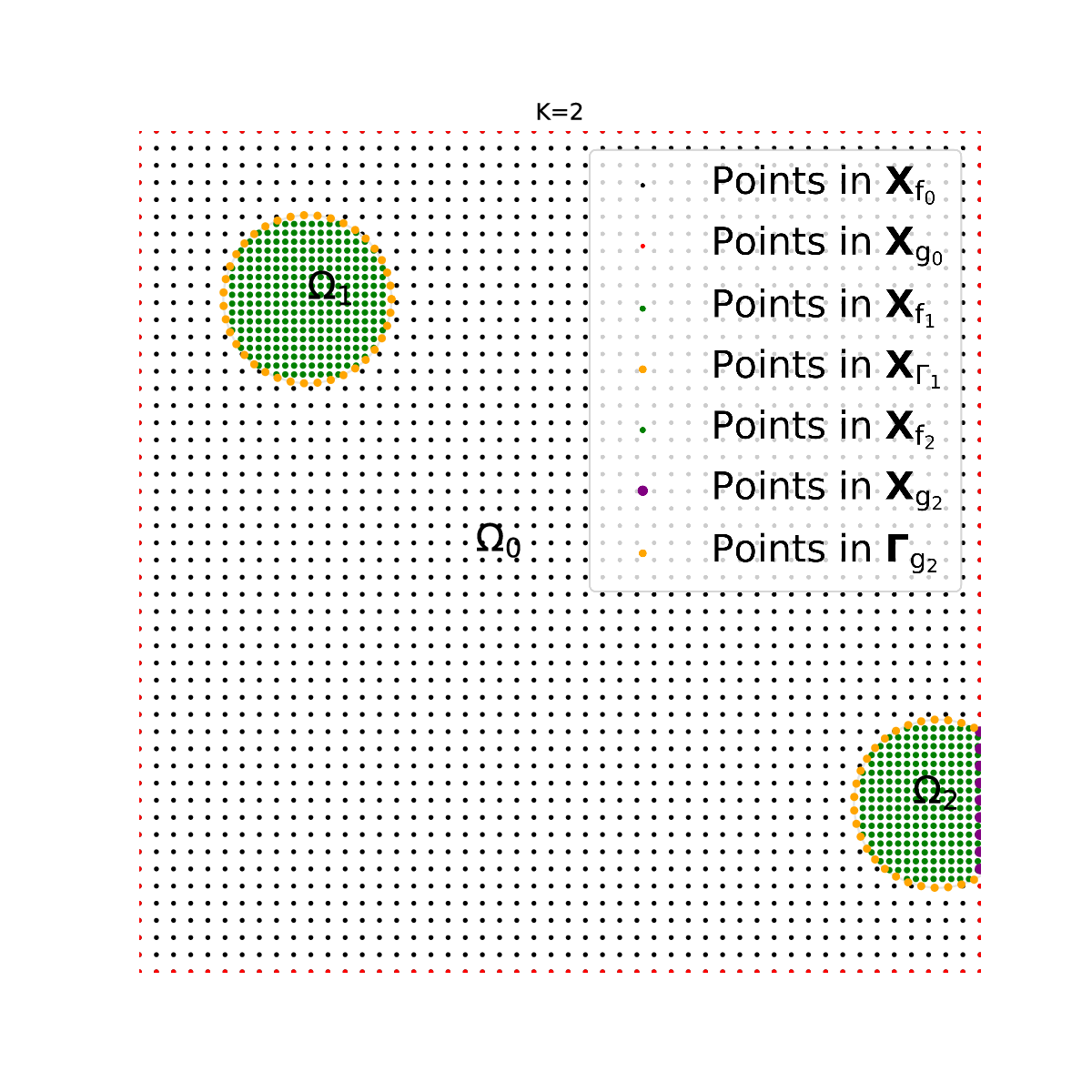}
\end{minipage}\\
\caption{A schematic diagram  from $K=0$ to $K=2$.}
\label{fig:K=2}
\end{figure}

In this section, we will elaborate on the procedure of the ANNB method.
Initially, we  denote $\Omega$ as $\Omega_{0}$. Assuming that $K~(K\ge 0)$  subdomains on which the solution of the PDE is low-regular have been identified and the domain $\Omega$ is partitioned into $K+1$ subdomains $\Omega_{0}, \Omega_{1}, \cdots, \Omega_{K}$, we denote by $\{\psi_{m,k}\}_{m=0}^{M_{k}}$ the neural network basis functions in subdomain $\Omega_{k}$. Additionally, we denote by $\bm{X}_{f_{k}},\bm{X}_{g_{k}}$ and $\bm{X}_{{\Gamma}_{k}}$ as the set of collocation points from $\Omega_{k}$, $\partial\Omega_{k}\cap\partial\Omega$ and $\Gamma_{k}$ respectively. For clarity, we illustrate the case of \(K=1\) within a two-dimensional square, as shown in the Figure \ref{fig:K=2}.

Then we get a least square problem \eqref{pro:MD-min} and the solution of this problem  can be obtained by applying Algorithm \ref{alg:DDM-nonlinear}. We denote the solution as
\begin{equation}\label{eq:sol-alpha}
\bm{\alpha}_{0}^{*} = (\alpha_{0,0}^{*}, \alpha_{1,0}^{*}, \cdots, \alpha_{M_{0},0}^{*})^{\transpose}, \quad \cdots, \quad \bm{\alpha}_{K}^{*} = (\alpha_{0,K}^{*}, \alpha_{1,K}^{*}, \cdots, \alpha_{M_{K},K}^{*})^{\transpose}.
\end{equation}

We consider that the solution of the target PDE in \(\Omega_{0}\) is low-regular near the locations where the residual
$| \mathcal{L} \tilde{u}_{0}(\bm{\alpha}_{0}^{*},\bm{x}_{f_{0}})  - f \left( \bm{x}_{f_{0}} \right)|
$
is large until the mean residual
$$
L_{\Omega_0}(\bm{\alpha}_{0}^{*}) =\frac{1}{|\bm{X}_{{f}_{0}}|} \sum_{\bm{x}_{f_{0}}\in\bm{X}_{f_{0}}} \left[ \mathcal{L} \tilde{u}_{0}(\bm{\alpha}_{0}^{*},\bm{x}_{f_{0}})  - f \left( \bm{x}_{f_{0}} \right) \right]^2.
$$
is smaller than a threshold $\epsilon$. Here, $|\bm{X}_{f_{0}}|$ is the number of collocation points in the interior of $\Omega_{0}$.

If $L_{\Omega_{0}}(\bm{\alpha}_{0}^{*}) > \epsilon$, we choose a sufficiently small \(r_ k> 0\), and determine a ball of radius \(r_k\) centered at
\begin{equation}\label{eq:ball-cen}
  \bm{x}_{K+1} = \argmax_{\bm{x}_{f_{0}} \in \bm{X}_{f_{0}}} \left|\mathcal{L} \left(\sum_{m=0}^{M_{0}} \alpha_{m,0}^{*} \psi_{m,0} \right) (\bm{x}_{f_{0}})- f(\bm{x}_{f_{0}}) \right|.
\end{equation}
We then partition the domain $\Omega_{0}$ into two subdomains $\Omega_{*}$ and $\Omega_{**}$ as described in Algorithm \ref{alg:Dissection-domain} and write $\Omega_{*}$  as $\Omega_{0}$ and $\Omega_{**}$ as $\Omega_{K+1}$. Moreover, we assume the subdomain $\Omega_{K+1}$ does not intersect with other subdomain $\Omega_{k}~(1\le k\le K)$ on which the solution of target PDE is low-regular. In the each subdomain \(\Omega_{k}~(1\le k\le K)\), the set of collation points and the neural network basis functions remain unchanged.
\begin{algorithm}[H]
    \caption{Domain Decomposition}
    \label{alg:Dissection-domain}
    \begin{algorithmic}[1]
     \STATE{ Input $\bm{\alpha}^{0,*}=(\alpha^{0,*}_{0},\alpha^{0,*}_{1},\cdots,\alpha^{0,*}_{M_{0}})^\transpose$, $\{\psi_{m}^{0}\}_{m=0}^{M_{0}} ,\bm{X}_{f_{0}},r,\Omega_{0}.$}
    \STATE{Compute $\bm{x}_{K+1}$ according to Eq.\eqref{eq:ball-cen}.

    Let $\Omega_{*}=\Omega\setminus B_{r}(\bm{x}_{K+1})$ and $\Omega_{**}=B_{r}(\bm{x}_{K+1})\cap \Omega$ where  $B_{r}(\bm{x}_{K+1}) : = \{ \bm{x},\|\bm{x}-\bm{x}_{K+1}\|_{2}<r\}$.}
    \RETURN $\Omega_{*},\Omega_{**}$
    \end{algorithmic}
\end{algorithm}
Next, we establish  the set of collocation points for the sub-domain \(\Omega_{K+1}\) and update the set of collocation points for \(\Omega_{0}\) as follows. For simplicity, we denote $K+1$ as $K$.

The set \(\bm{X}_{g_{K}}\) consists of points from \(\bm{X}_{g_{0}}\) that are located on the boundary \(\partial \Omega_{K}\). For clarity, \(\bm{X}_{f_0}\) represents the collocation points from \(\bm{X}_{f_0}\) that fall within the domain \(\Omega_{0}\), and \(\bm{X}_{g_0}\) denotes the collocation points from \(\bm{X}_{g_0}\) that are located on the boundary \(\partial \Omega_{0}\). This process can be summarized as
$$
\bm{X}_{g_{k}} = \bm{X}_{g_{0}} \cap \partial \Omega_{K}, \quad \bm{X}_{f_0} = \bm{X}_{f_{0}} \setminus \Omega_{K}, \quad \bm{X}_{g_0} = \bm{X}_{g_{0}} \setminus \partial \Omega_{K}.
$$
Additionally, $\bm{X}_{f_{k}}$  consist of uniformly grid  points within the domain $\Omega_{K}$ and  $\bm{X}_{\Gamma_{K}}$ consist of the uniform distributed points on the interface $\Gamma_{K}$, respectively.

We next show how to generate neural network basis functions in the subdomain $\Omega_{K}= B_{r_{K}}(\bm{x}_{K})$ $\cap \Omega$. For a positive integer $M_{k}$, based on the standard neural network basis functions
$$
\psi_{0,K}(\bm{x}) = 1, \quad \psi_{m,K}(\bm{x}) = \sigma({\bm{\omega}_{m,K}}^{\transpose}\bm{x} + b_{m,K}) \quad (1 \le m \le M_{k}),
$$
we form a scaled neural network basis functions as follows.
\begin{equation}\label{eq:MSNN}
\psi_{0,K}(\bm{x}) = 1, \quad \psi_{m,K}(\bm{x},c_{K}) = \sigma\left(c_{K} \cdot {\bm{w}_{m,K}}^{\transpose} (\bm{x} - \bm{x}_{K}) + b_{m,K}\right) \quad (1 \le m \le M_{k}),
\end{equation}
where $c_{K}$ is the scaling coefficient to be determined below.

Next, we represent \( u_{K}(\bm{x}) \) as a function spanned by the neural network basis functions, as shown in Eq. \eqref{eq:MSNN}. To determine the scaling coefficient \( c_{K} \), we formulate the following least square problem based on the strong formulation of the last four equations in \eqref{pro:MD} at the collocation points \( \bm{X}_{f_K} \), \( \bm{X}_{g_{K}} \), and \( \bm{X}_{\Gamma_{K}} \):

\begin{equation}\label{pro:K-min}
\begin{aligned}
  \argmin_{\alpha_{K}c_{K}}L_{K}(\bm{\alpha}_{K},c_{K})=\min_{\bm{\alpha}_{K},c_{K}} &\left\{
 \sum_{\bm{x}_{f_{K}}\in \bm{X}_{f_{k}}} \left( \mathcal{L}\left( \sum_{m=0}^{M_{K}} \alpha_{m,K} \psi_{m,K}\right)(\bm{x}_{f_{K}},c_{K})  - f\left(\bm{x}_{f_{K}}\right) \right)^2  \right.\\
& +\sum_{\bm{x}_{g_{K}}\in\bm{X}_{g_{k}}}\left( \sum_{m=0}^{M_{K}} \alpha_{,K} \psi_{m,k}(\bm{x}_{g_{K}},c_K) -g(\bm{x}_{g_K})\right)^2 \\
&+ \sum_{\bm{x}_{\Gamma_{K}}\in\bm{X}_{\Gamma_{K}}}\left[ \left( \sum_{m=0}^{M_{K}} \alpha_{m,K} \psi_{m,K}(\bm{x}_{\Gamma_{K}},c_{K})- \sum_{m=0}^{M_{0}} \alpha_{m,0}^{*} \psi_{m,0}(\bm{x}_{\Gamma_{K}})  \right)^2\right.\\
&\left.\left.+ \left( \sum_{m=0}^{M_{K}} \alpha_{m,K} \frac{\partial\psi_{m,K}(\bm{x}_{\Gamma_{K}},c_{K})}{\partial\bm{n}_{K}} - \sum_{m=0}^{M_{0}} \alpha_{m,0}^{*} \frac{\partial\psi_{m,0}(\bm{x}_{\Gamma_{K}})}{\partial\bm{n}_{K}} \right)^2\right]\right\}.
\end{aligned}
\end{equation}

The solution $c_{K}^{*}$ of problem \eqref{pro:K-min} is selected as the scaling coefficient for $\Omega_{K}$.

To simplify the computational cost, we choose the scaling coefficient $c_K$ from an integer between 1 and an positive integer $L$ (say, $L=10$), cf. Algorithm \ref{alg:scale-basis} for further details.  Moreover, for the help of readers, we carry out the solution procedure of problem \eqref{pro:min-scale factor-s} in the appendix. Additionally, the dichotomy method can also be employed to determine the scaling coefficients.
\begin{algorithm}[H]
  \caption{Generate neural basis functions for subdomain $\Omega_{K}$.}
  \label{alg:scale-basis}
  \begin{algorithmic}[1]
  \STATE{$\mathbf{Input}$ subdomain $\Omega_{K}=\mathrm{B}_{r_{K}}(\bm{x}_{c_{K}})\cap \Omega$; number of neural network basis $M_{K}$;~
   an approximate solution $\tilde{u}_{0}(\bm{x})=\sum_{m=0}^{M_{0}}\alpha_{m,0}^{*}\psi_{m,0}(\bm{x})$;~maximum number of iteration steps $N$;~tolerance $tol$;~$n=1$.}
   \STATE{$\mathbf{Output}:
   \{\psi_{m,K}(\bm{x})\}_{m=0}^{M_K}$ where
   $\psi_{0,K}=1,~
\psi_{M_{K}}(\bm{x})=\sigma\left(s^{*}*\bm{w}_{m,K}^{\top} (\bm{x}-\bm{x}_{{K}})+b_{m,K}\right)~(1\le m\le M_{K})$.}
   \FOR{s=1,2,3,4,5,6,7,8,9,10}
  \FOR{$m=1:M_{K}$}
  \STATE{
Generate weights and biases: $\bm{w}_{m,K}, b_{m,K}$.
\\
  Construct the neural basis function: $\psi_{M_{K}}(\bm{x})=\sigma\left(s*{\bm{w_{m,K}}^{\top}} (\bm{x}-\bm{x}_{{K}})+b_{m,K}\right)$}

  \ENDFOR
\STATE{
Denote $\bm{\alpha}_{K}^{*}$ as  the solution of problem
\begin{equation}\label{pro:min-scale factor-s}
\argmin_{\alpha_{K}}L_{K}(\bm{\alpha}_{K},s).
\end{equation}
The details of solving problem \eqref{pro:min-scale factor-s} are shown in Appendix.}

\STATE{Let $\Loss_{s}=L_{K}(\bm{\alpha}_{K}^{*},s)
.$}
\ENDFOR
\STATE{Let $s^*= \argmin_{1\le s\le 10}\Loss_{s}$.}
  \end{algorithmic}
\end{algorithm}

With the basis function on each subdomain $\Omega_k$ determined, we finally apply  Algorithm \ref{alg:DDM-nonlinear} to solve the new least square problem \eqref{pro:MD-min} and obtain new solution
\begin{equation*}
\bm{\alpha}_{0}^{*} = (\alpha_{0,0}^{*}, \alpha_{1,0}^{*}, \cdots, \alpha_{,0M_{0}}^{*})^{\transpose}, \quad \cdots, \quad \bm{\alpha}_{K}^{*} = (\alpha_{0,K}^{*}, \alpha_{1,K}^{*}, \cdots, \alpha_{M_{k},K}^{*})^{\transpose}.
\end{equation*}

This iterative partitioning of the domain $\Omega_{0}$ continues until $L_{\Omega_{0}}(\bm{\alpha}_{0}^{*}) < \epsilon$. Then the approximate solution of the original PDE can be described as follows.
\begin{equation}\label{eq:app-u}
u(x)=\Tilde{u}_{k}(\bm{\alpha}_{k}^{*},\bm{x}) = \sum_{m=0}^{M_{k}} \alpha_{m,k}^{*} \psi_{m,k}(\bm{x})~ \bm{x}\in \Omega_{k},\quad 0\le k\le K
\end{equation}
The main steps of the ANNB method are summarized in Algorithm \ref{alg:ANNB}.

\begin{algorithm}[H]
    \caption{ANNB method for second-order PDE with low-regular solutions.}
    \label{alg:ANNB}
    \begin{algorithmic}[1]
    \STATE{Let $~K=0,~\Omega_{0}=\Omega.$}
    \STATE{$\mathbf{Input}$:~$M_{0}$, number of neural network basis functions,;~$\bm{X}_{f_{0}}$~a set of collocation points in $\Omega_{0}$,;$\bm{X}_{g_{0}}$~a set of collocation points on $\partial\Omega_{0}\cap\partial\Omega_{0}$;~$\epsilon$.}
     \STATE{$\mathbf{Output}: u(\bm{x})$.}
     \FOR{$m=1:M_{0}$}
    \STATE{
    Generate weight and bias: $\bm{w}_{m,0}, b_{m,0}$.

    Construct the neural basis function: $\psi_{m,0}(\bm{x})=\sigma\left({\bm{w}_{m,0}}^{\top} \bm{x}+b_{m,0}\right).$
}
\ENDFOR
    \vspace{0.1cm}
    \STATE{  Obtain $\bm{\alpha}_{k}^{*}=(\alpha_{0,k}^{*},\alpha_{1,k}^{*},\cdots,\alpha_{M_{k},k}^{*})^{\transpose}~(0\le k\le K)$ obtained by applying  Algorithm\ref{alg:DDM-nonlinear}  for solving the minimization problem \eqref{pro:MD-min}.}

    \WHILE{$L_{\Omega_{0}}(\bm{\alpha}_{0}^{*})>\epsilon$}
        \vspace{0.1cm}
    \STATE{
      Partition $\Omega_{0}$ into $\Omega_{*}$ and $\Omega_{**}$ by applying Algorithm\ref{alg:Dissection-domain}.}
        \vspace{0.1cm}
    \STATE{Let $K=K+1;~\Omega_{0}=\Omega_{*};~\Omega_{K}=\Omega_{**};~\bm{X}_{f_{0}} = \bm{X}_{f_{0}}\cap \Omega_{0};~\bm{X}_{g_{K}}=\bm{X}_{g_{0}}\cap\partial\Omega_{K};~\bm{X}_{g_{0}}=\bm{X}_{g_{0}}\cap\partial\Omega_{0}.$  And generate $\bm{X}_{f_{K}},~\bm{X}_{\Gamma_{K}}.$}

    \vspace{0.1cm}
\STATE{ Obtain  $\{\psi_{M_{K}}(\bm{x})\}_{m=0}^{M_{K}}$ by applying Algorithm \ref{alg:scale-basis}.}

    \vspace{0.1cm}
\STATE{Obtain $\bm{\alpha}_{k}^{*}=(\alpha_{0,k}^{*},\alpha_{1,k}^{*},\cdots,\alpha_{M_{k},k}^{*})^{\transpose}~(0\le k\le K)$ obtained by applying  Algorithm\ref{alg:DDM-nonlinear}  for solving the minimization problem \eqref{pro:MD-min}.}
\ENDWHILE
\STATE{Let
\begin{equation}
u(\bm{x})\approx\Tilde{u}(\bm{\alpha}_{k}^{*},\bm{x})=\sum_{m=0}^{M_{k}} \alpha_{m,k}^{*} \psi_{m,k}(\bm{x}) \quad \bm{x}\in\Omega_{k}.
\end{equation}}
\end{algorithmic}
\end{algorithm}
\begin{remark}
       It is worth noting that Algorithm \ref{alg:ANNB}  is not only applicable to neural network basis functions but can applies  to any class of mesh-free basis functions (cf. \cite{2008-meshless}) as well.
\end{remark}
For clarity, we present an example from \( K=0 \) to \( K=2 \) in Fig. \ref{fig:K=2} to show the process of generating the collocation points. In the left picture, \(\Omega = \Omega_{0} = [0,1]^2\), the set \( \bm{X}_{f_{0}} \) consists of black points, and the set \( \bm{X}_{g_{0}} \) consists of red points. In the middle picture, \( \Omega \) is partitioned into two subdomains, $\Omega_{0}$ and $\Omega_{1}$, where \( \Omega_{1} = B_{0.1}(\bm{x}_{1}) \), with \( \bm{x}_{1} = [0.2,0.8]^\transpose \), and \( \Omega_{0} = \Omega \setminus \Omega_{1} \). The set \( \bm{X}_{f_{1}} \) consists of green points, the set \( \bm{X}_{\Gamma_{1}} \) consists of orange points, and the set \( \bm{X}_{g_{1}} \) consists of the points from \( \bm{X}_{g_{0}} \) that fall on the boundary \( \partial\Omega_{1}\),~i.e., $\bm{X}_{g_{1}}=\bm{X}_{g_{0}}\cap\partial\Omega_{1}=\emptyset$ . Furthermore, the set \( \bm{X}_{f_{0}} \) consists of the points from \( \bm{X}_{f_{0}} \) that fall within the subdomain \( \Omega_{0} \),~i.e.,~\( \bm{X}_{f_{0}} = \bm{X}_{f_{0}} \cap \Omega_{0} \), and the set \( \bm{X}_{g_{0}} \) consists of the points from \( \bm{X}_{g_{0}} \) that fall on the boundary \( \partial\Omega_{0} \),~i.e.,~\( \bm{X}_{g_{0}} =\bm{X}_{g_{0}}\cap\partial\Omega_{0} \).

In the right picture, \( \Omega \) is partitioned into three subdomains, $\Omega_{0}$, $\Omega_{1}$ and $\Omega_{2}$, where \( \Omega_{2} = B_{0.1}(\bm{x}_{2})\) with \( \bm{x}_{2} = [0.95,0.2]^\transpose \)and \( \Omega_{0} = \Omega_{0} \setminus \Omega_{2} \). The set $\bm{X}_{f_{1}},\bm{X}_{g_{1}}$ and $\bm{X}_{\Gamma_{1}}$ remain unchanged. The set \( \bm{X}_{f_{2}} \) consists of the green points in \( \Omega_{2} \), the set \( \bm{X}_{\Gamma_{2}} \) consists of the orange points on the boundary of \( \Omega_{2} \) and the set $\bm{X}_{g_{2}}$ consists of the purple points on the boundary of \( \Omega_{2} \), i.e., \(\bm{X}_{g_{2}}=\bm{X}_{g_{0}}\cap\partial\Omega_{2}\). Moreover,~\( \bm{X}_{f_{0}} = \bm{X}_{f_{0}} \cap \Omega_{0} \) and \( \bm{X}_{g_{0}} =\bm{X}_{g_{0}}\cap\partial\Omega_{0} \).

\section{Numerical experiments}
In all numerical experiments, we generate the neural network basis functions  using the second strategy in Section \ref{sec:TNNs}.

\subsection{Experiment Setup}
Given the test set \( S_{\mathrm{test}} = \{\bm{x}_{i}\}_{i=1}^{N_{t}} \), we evaluate the performance of our algorithms using the following error metric:

\begin{equation}\label{eq:re_l2err}
\operatorname{err}_{L_2} = \frac{\sqrt{\sum_{i=1}^{N_{t}} \left|\hat{u}\left(\bm{x}_i\right) - u\left(\bm{x}_i\right)\right|^2}}{\sqrt{\sum_{i=1}^{N_{t}} \left|u\left(\bm{x}_i\right)\right|^2}},
\end{equation}
where \( N_{t} \) is the total number of test points, \( \hat{u}\left(\bm{x}_i\right) \) denotes the predicted solution value, and \( u\left(\bm{x}_i\right) \) denotes the exact solution value at test point \( \bm{x}_i \).


Here is how we construct \( \bm{X}_{f_{k}} \) and \( \bm{X}_{\Gamma_{k}} \) (where \( k \geq 1 \)). We define \( \bm{x}_{k} \) as \( (x_{c_{k}}, y_{c_{k}})^{\transpose} \) for all two-dimensional examples, and as \( (x_{c_{k}}, y_{c_{k}}, z_{c_{k}})^{\transpose} \) for three-dimensional examples.

\begin{itemize}
    \item In all two-dimensional examples, \( \bm{X}_{f_{k}} \) consists of a \( 40 \times 40 \) uniform grid points on \( [x_{c_{k}} - r_{k}, x_{c_{k}} + r_{k}] \times [y_{c_{k}} - r_{k}, y_{c_{k}} + r_{k}] \), with grid points outside the sub-domain \( \Omega_{k} \) masked off. \( \bm{X}_{\Gamma_{k}} \) includes 200 uniformly distributed points on \( \partial B_{r_{k}}(\bm{x}_{c_{k}}) \), with points not on the interface \( \Gamma_{k} \) masked off.
    \item In the three-dimensional example \ref{sec:ex-3d}, \( \bm{X}_{f_{k}} \) is constructed with 8500 points uniformly distributed within \( [x_c - r_{k}, x_c + r_{k}] \times [y_c - r_{k}, y_c + r_{k}] \times [z_c - r_{k}, z_c + r_{k}] \), with points outside \( \Omega_{k} \) masked off. \( \bm{X}_{\Gamma_{k}} \) includes 600 points uniformly distributed on \( \partial B_{r_{k}}(\bm{x}_{c_{k}}) \), with points not on the interface \( \Gamma_{k} \) masked off.
\end{itemize}
Additionally, in algorithms \ref{alg:DDM-nonlinear} and \ref{alg:scale-basis}, we set \( tol = 1e-5 \). In all examples, the number of basis functions in each subdomain $\Omega_{k}~(k\ge 1)$ is set to \( M^{*} \).

\subsection{Two-dimensional problem with one peak and multiple peaks.}\label{sec:linear-2d-peak}
We first consider the following two-dimensional Poisson equation \cite{2023-FIPINN-zhou1,2023-FIPINN-zhou2}:
\begin{equation}\label{eq:poisson2d}
 \begin{aligned}
-\Delta u(x, y)=f(x, y)  &&  \text { in } \Omega, \\
u(x, y)=g(x, y)          && \text { on } \partial \Omega,
\end{aligned}
\end{equation}
where $\Omega$ is $[-1,1]^2$ and we specify the true solution as
$$
u(x, y)=\sum_{p=1}^{P}\exp \left(-1000\left[(x-x_{p})^2+(y-y_{p})^2\right]\right),
$$
which has several peaks at $\left(x_p, y_p\right), p=1, \ldots, P$, and will decay to zero exponentially in other places. We considered three cases in this part:
\begin{enumerate}
\item Case 1: One peak located at $ (0.5, 0.5) $.
\item Case 2: Peaks located at $ (0.5, 0.5) $ and $ (-0.5, -0.5) $.
\item Case 3: Peaks located at $ (\pm 0.5, \pm 0.5) $.
\end{enumerate}

In algorithm \ref{alg:ANNB}, we set $\epsilon =1e-04,~M_{0}=200$ and $r_{k}=r=0.15~(k\ge1)$.  The data set $\bm{X}_{f_{0}}$ consists of the grid points within a $50 \times 50$ uniform grid points in $[-1,1]^2,$ i.e., $J_{f_0}=2500$. The data set $\bm{X}_{g_{0}}$ consists of 400 uniformly distributed points on $\partial \Omega ,$ i.e., $J_{g_{0}}=400.$ Additionally, we compute the $ \operatorname{err}_{L_{2}} $ \eqref{eq:re_l2err} on a $ 256 \times 256 $ uniform gird points in $ \Omega = [-1,1]^2 $.

\begin{figure}[H]
  \begin{minipage}{0.3\textwidth}
    \centering
    \includegraphics[width=\textwidth, height=5cm]{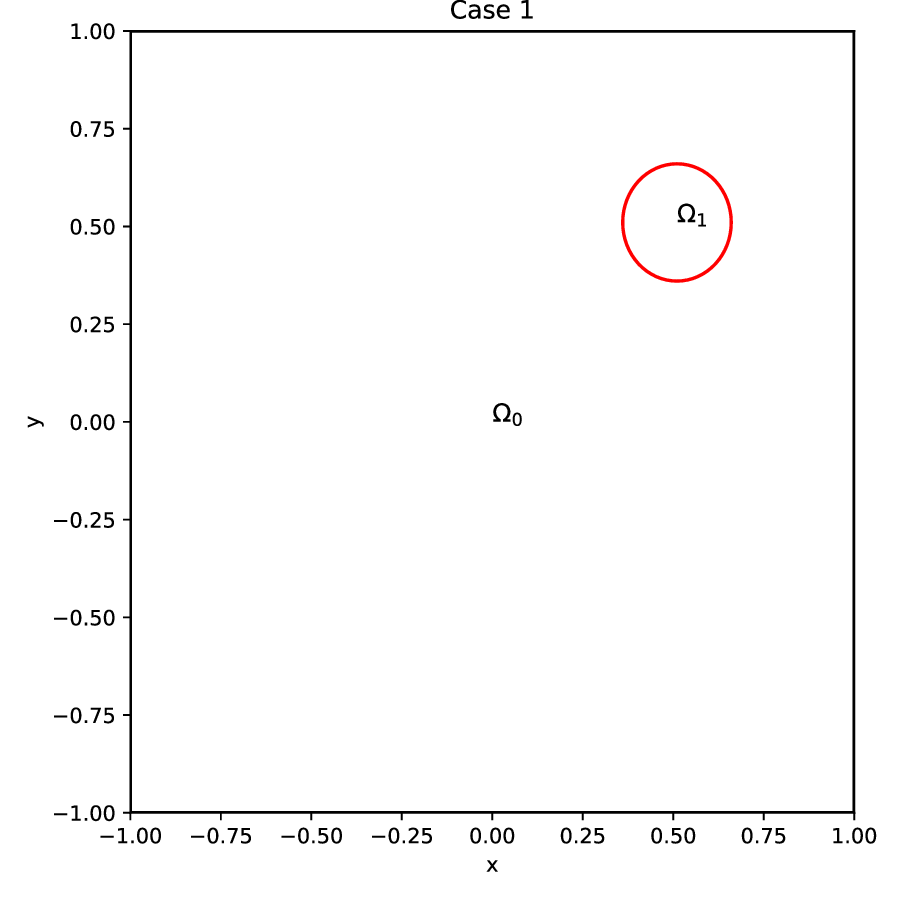}
  \end{minipage}
  \begin{minipage}{0.3\textwidth}
    \centering
    \includegraphics[width=\textwidth, height=5cm]{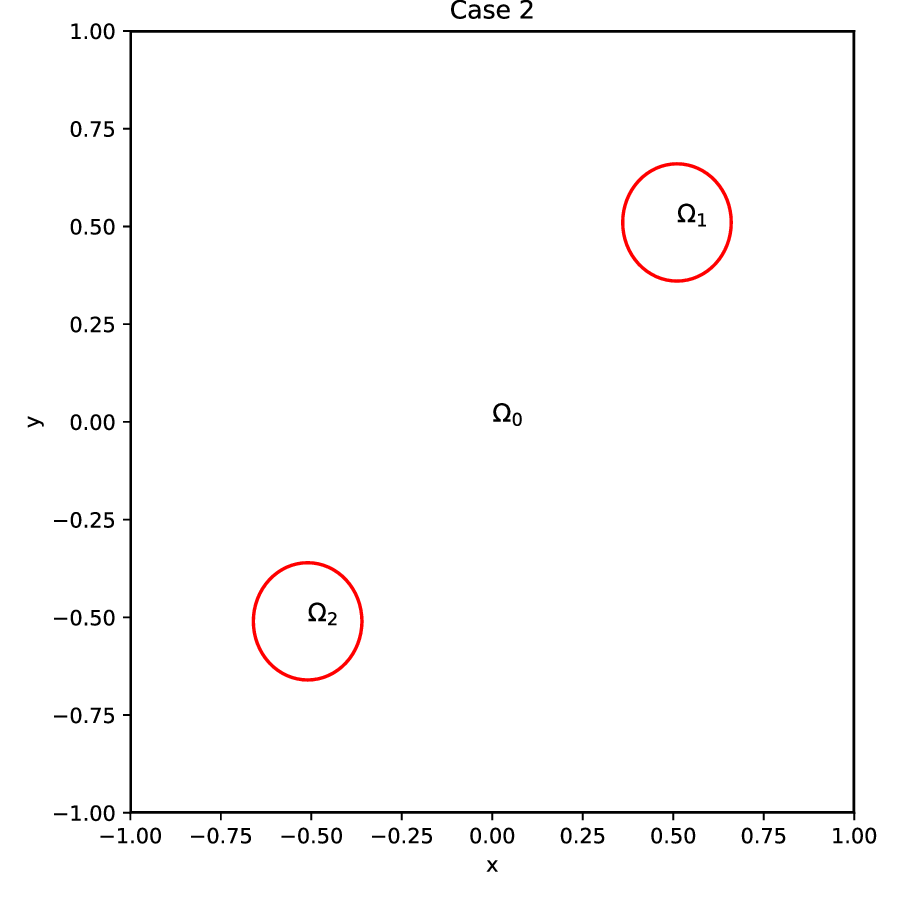}
  \end{minipage}
  \begin{minipage}{0.3\textwidth}
    \centering
    \includegraphics[width=\textwidth, height=5cm]{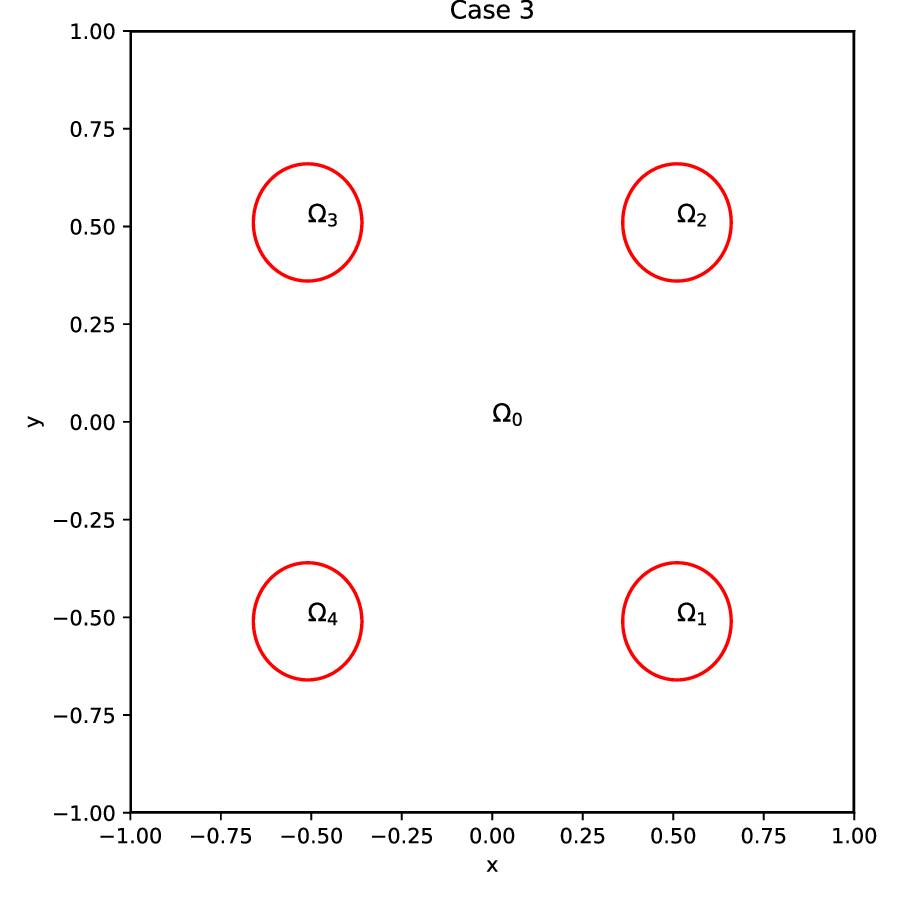}
  \end{minipage}
  \caption{The schematic diagrams of the domain for example \ref{sec:linear-2d-peak}.}
  \label{Fig:Domain}
\end{figure}

\begin{figure}[H]
  \begin{minipage}{0.33\textwidth}
    \centering
    \includegraphics[width=\textwidth, height=5cm]{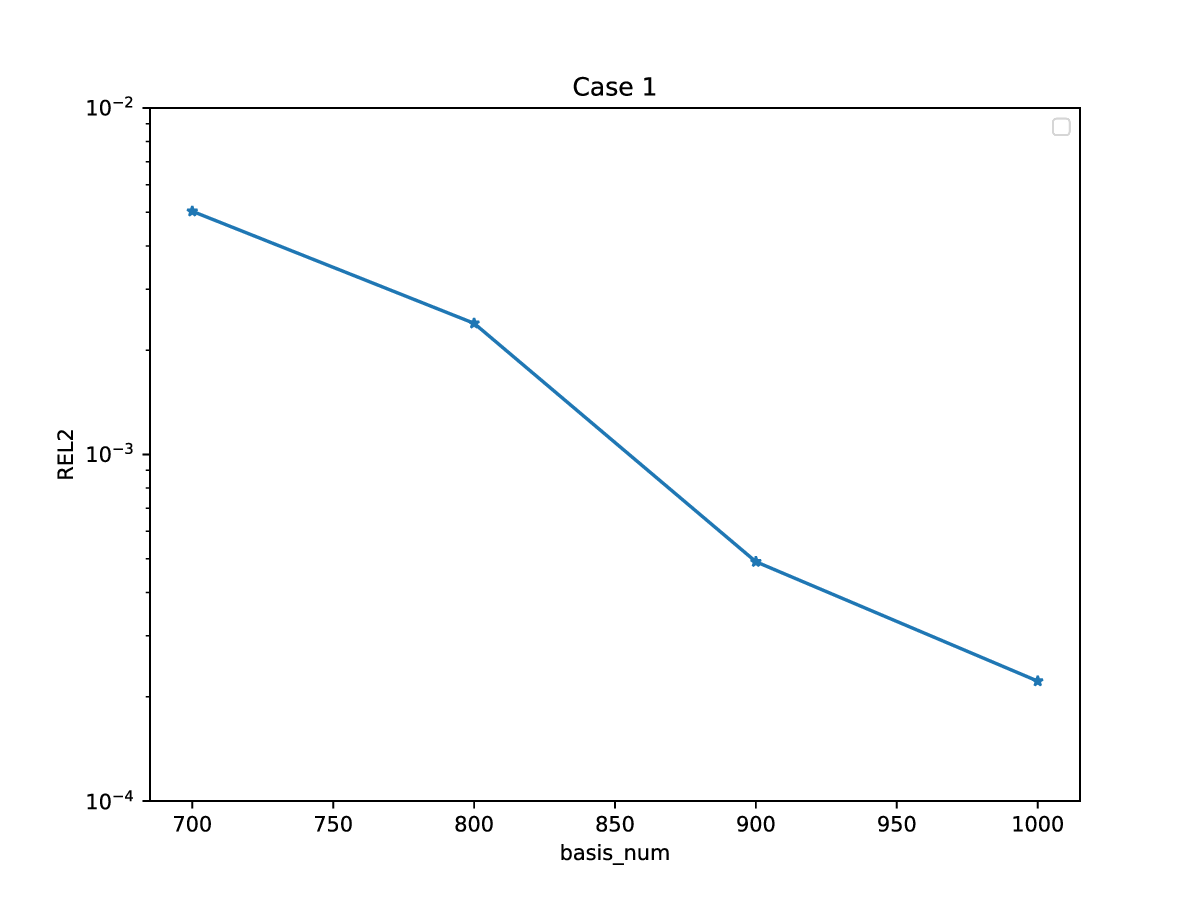}
  \end{minipage}
  \begin{minipage}{0.33\textwidth}
    \centering
    \includegraphics[width=\textwidth, height=5cm]{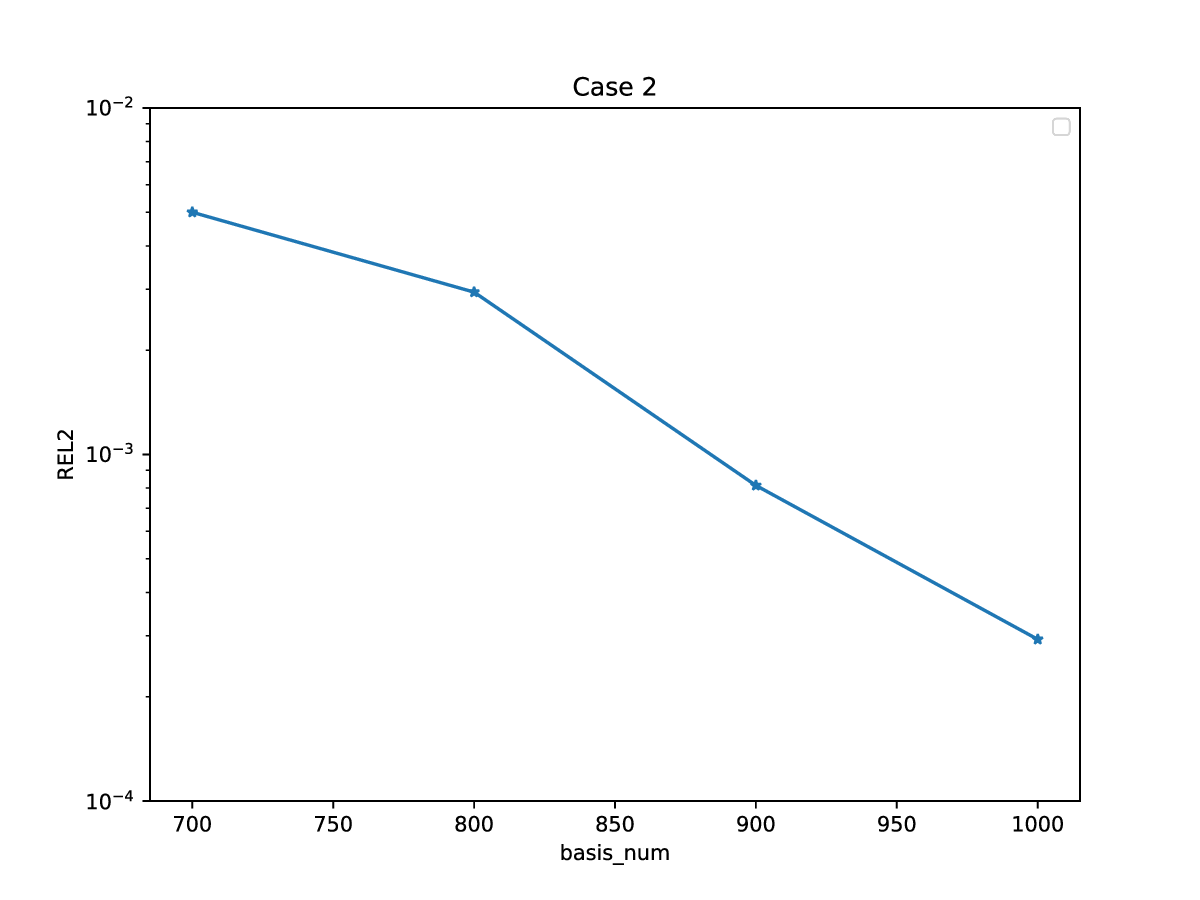}
  \end{minipage}
  \begin{minipage}{0.33\textwidth}
    \centering
    \includegraphics[width=\textwidth, height=5cm]{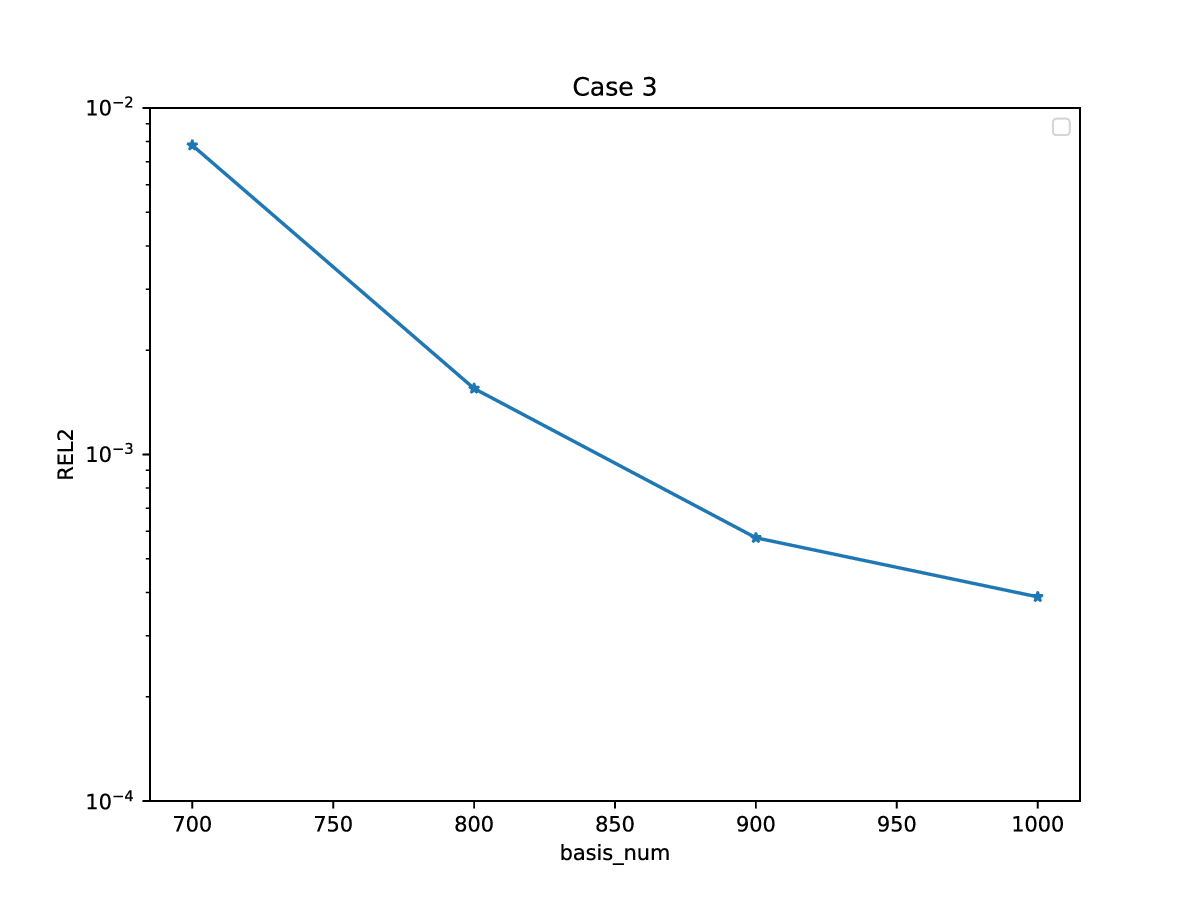}
  \end{minipage}\\
  \caption{The $\operatorname{err}_{L_{2}}$ with the number of scaled neural network basis functions.}
  \label{Fig:rel2}
\end{figure}

\begin{figure}[H]
  \begin{minipage}{\textwidth}
    \centering
    \includegraphics[width=\linewidth,height=5cm]{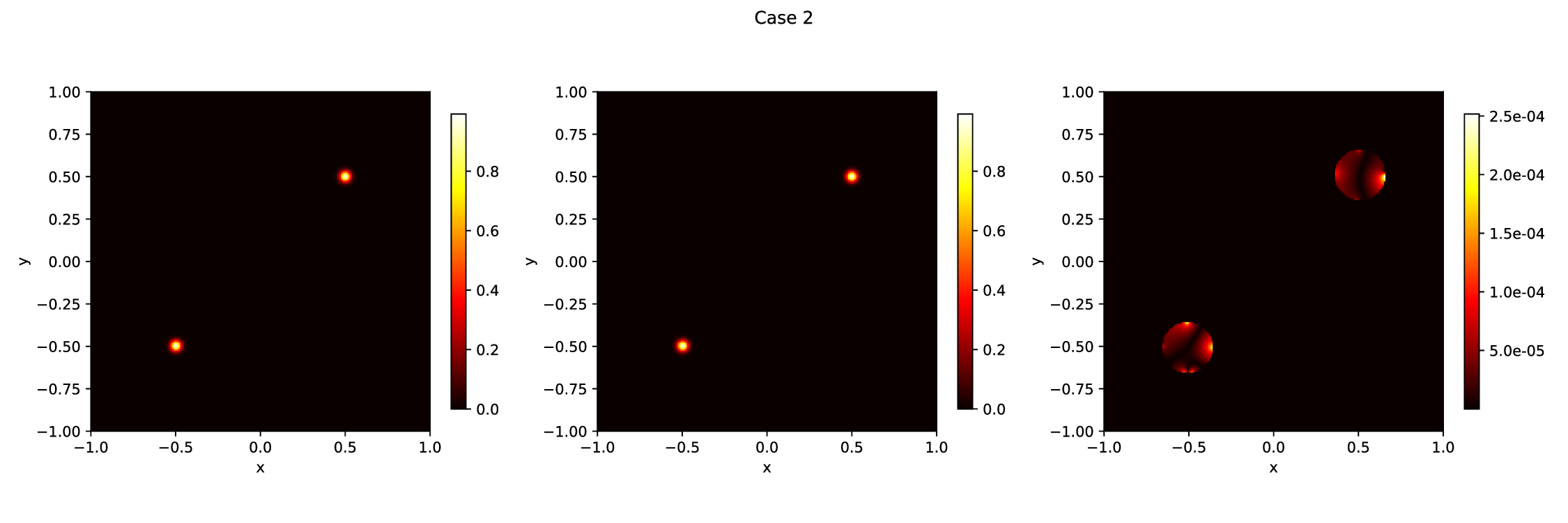}
  \end{minipage}\\
  \begin{minipage}{\textwidth}
    \centering
    \includegraphics[width=\linewidth,height=5cm]{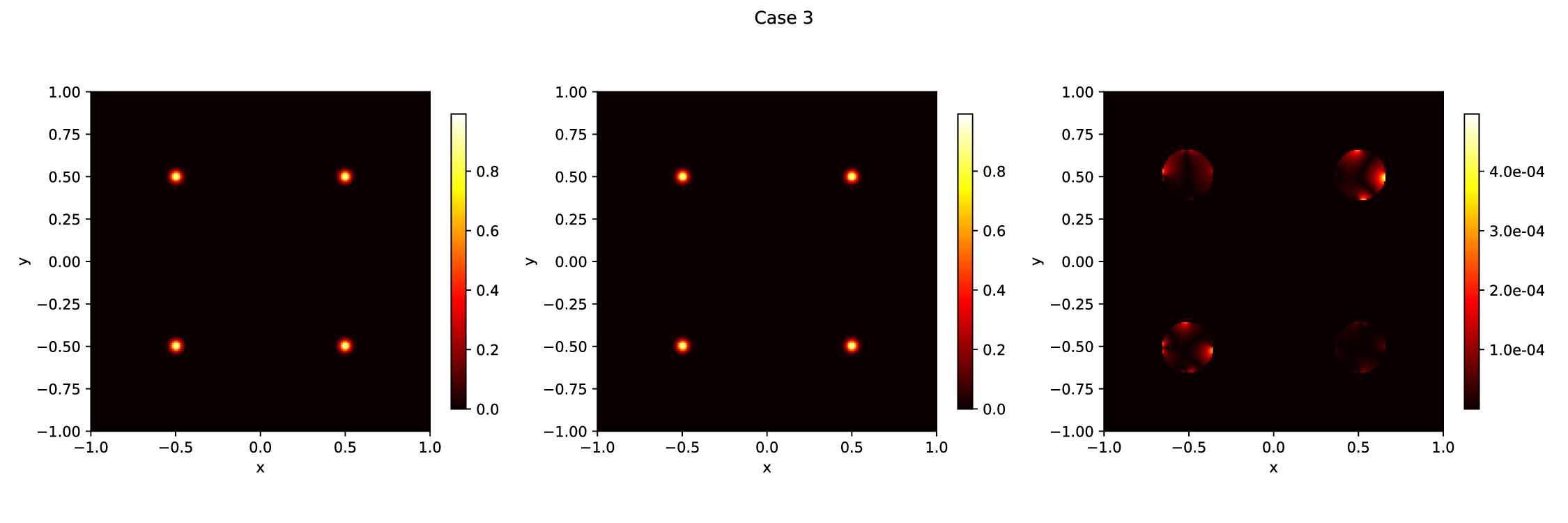}
  \end{minipage}\\
  \caption{The exact solution~(Left), the predicted solution~(Middle) and absolute error~(Right) when $M^{*}=1000$.}
  \label{Fig:heat_cap}
\end{figure}

The numerical results are shown in Figs. \ref{Fig:Domain}, \ref{Fig:rel2}, \ref{Fig:heat_cap}  and Tab. \ref{tab:scale-facotr-ex1},  respectively.
 We can clearly see from Fig. \ref{Fig:Domain} that
the domain $\Omega=[-1,1]^2$ is eventually partitioned into two sub-domains for case 1, three sub-domains for case 2 and five sub-domains for case 3.  The following is the information about sub-domains with low-regular solutions  for three cases.
\begin{itemize}
    \item Case1: $\Omega_{1}=B_{r}(\bm{x}_{1})$, where $\bm{x}_{1}=(0.5102,0.5102)^{\transpose}$.
    \item Case2: $\Omega_{k} = B_{r}(\bm{x}_{k}), k=1, 2$, where $\bm{x}_{1}=(-0.5102,-0.5102)^{\transpose},~\bm{x}_{2}=(0.5102,0.5102)^{\transpose}$.
    \item Case3: $\Omega_{k}=B_{r}(\bm{x}_{k}),k=1, 2, 3, 4$, where
   $\bm{x}_{1}=(-0.5102,0.5102)^{\transpose}, \bm{x}_{2}=(-0.5102,-0.5102)^{\transpose},~\bm{x}_{3}=(0.5102,-0.5102)^{\transpose},~\bm{x}_{4}=(0.5102,0.5102)^{\transpose}$.
\end{itemize}

For three cases, it is clear that  the number of  sub-domains with low-regular solutions is equal to the number of peaks, and the positions of these sub-domains coincide with the positions of the peaks. Besides, we observe from Tab. \ref{tab:scale-facotr-ex1} that the scaling coefficient on each subdomain with low-regular solution obtained by algorithm $\ref{alg:scale-basis}$  equal 5 for three cases.


\begin{table}[H]
    \centering
    \begin{tabular}{c|c|c|c}
                & Case 1 & Case 2 & Case 3  \\\hline
        $c_{1}$& 5        & 5        &5 \\\hline
        $c_{2}$& $\times$ & 5        &5 \\\hline
        $c_{3}$& $\times$ & $\times$ &5 \\\hline
        $c_{4}$& $\times$ &$\times$  &5 \\\hline
    \end{tabular}
    \caption{Scale coefficients obtained by Algorithm $\ref{alg:scale-basis}$ for three cases. }
    \label{tab:scale-facotr-ex1}
\end{table}

Fig. \ref{Fig:rel2} shows, in three cases, the trend of $\operatorname{err}_{L_{2}}$ defined in \eqref{eq:re_l2err} decrease with the increasing number of neural network basis functions $M^{*}$, where $M^{*}$ is chosen as $700, 800, 900$ and $1000$. We observe that, in three cases, our method exhibits superior performance, and the $\operatorname{err}_{L_{2}}$ obtained by ANNB reaches approximately  $1\mathrm{e}-4$ when $M^*=1000$, which is significantly smaller than the results in\cite{2023-FIPINN-zhou1,2023-FIPINN-zhou2}.

Fig. \ref{Fig:heat_cap} displays the predicted solutions and their corresponding absolute error for the three cases when $M^{*} = 1000 $. From Fig. \ref{Fig:heat_cap}, it is evident that the value of absolute error within each sub-domain $ \Omega_{k}$ ($ 1 \leq k \leq 4 $) are relatively small, although it is larger than that in the domain $\Omega_{0} $.

\begin{figure}[H]
  \begin{minipage}{\textwidth}
    \centering
    \includegraphics[width=\linewidth,height=5cm]{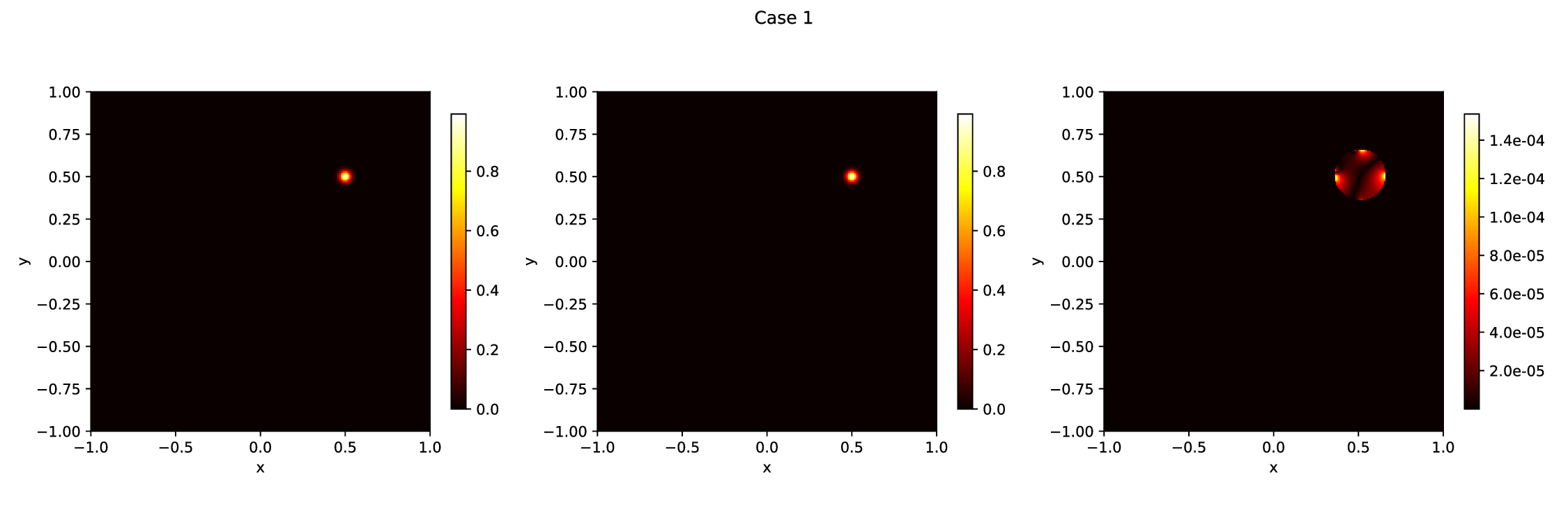}
  \end{minipage}\\
\end{figure}

\subsection{Two dimensional nonlinear problem with one peak and multiple peaks}\label{sec:nonlinear-peaks}
Consider the following two-dimensional nonlinear problem:
\begin{equation}\label{eq:nonlinear-poisson2d}
 \begin{aligned}
-\Delta u(x, y) +u(x,y)^{2}=f(x, y) & & \text { in } \Omega, \\
u(x, y)=g(x, y) & & \text { on } \partial \Omega,
\end{aligned}
\end{equation}
where $\Omega$ is $[-1,1]^2$ and we specify the true solution as
$$
u(x, y)=\sum_{p=1}^{P}\exp \left(-1000\left[(x-x_{p})^2+(y-y_{p})^2\right]\right).
$$
In this example,  we also consider the same three cases as in example \ref{sec:linear-2d-peak}. In addition, we also calculate the $\operatorname{err}_{L_{2}}$  \eqref{eq:re_l2err} on a $256 \times 256$ uniform grid points in $\Omega=[-1,1]^2$. In algorithm \ref{alg:ANNB}, we set $\epsilon =1\mathrm{e}-4,~M_{0}=200$
and $r_{k}=r=0.15~(k\ge 1).$  The data set $\bm{X}_{f_{0}}$ and $\bm{X}_{g_{0}}$ are same as the example \eqref{sec:linear-2d-peak}.

\begin{figure}[H]
  \begin{minipage}{0.3\textwidth}
    \centering
    \includegraphics[width=\textwidth, height=5cm]{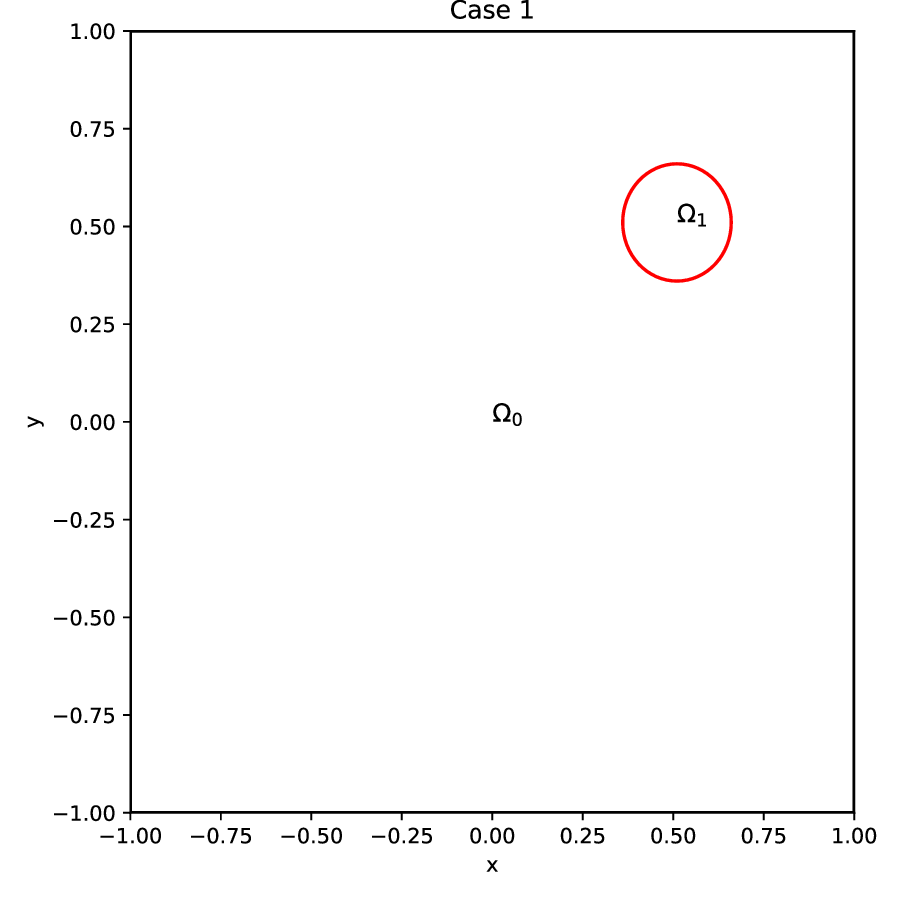}
  \end{minipage}
  \begin{minipage}{0.3\textwidth}
    \centering
    \includegraphics[width=\textwidth, height=5cm]{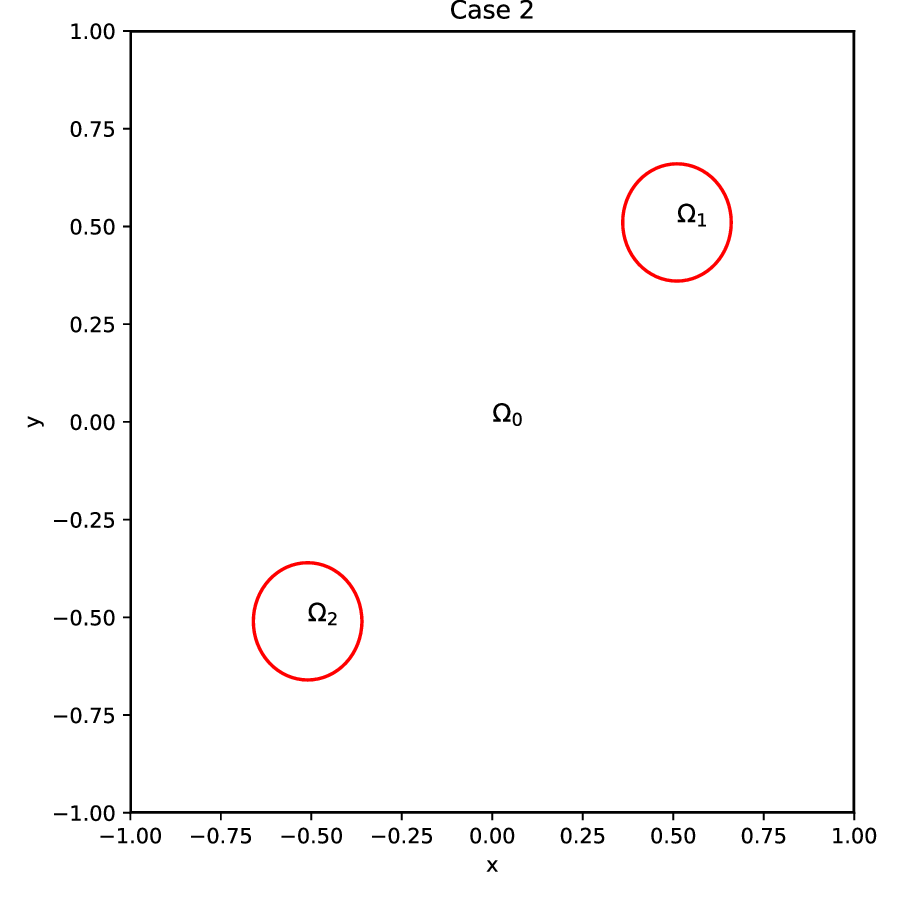}
  \end{minipage}
  \begin{minipage}{0.3\textwidth}
    \centering
    \includegraphics[width=\textwidth, height=5cm]{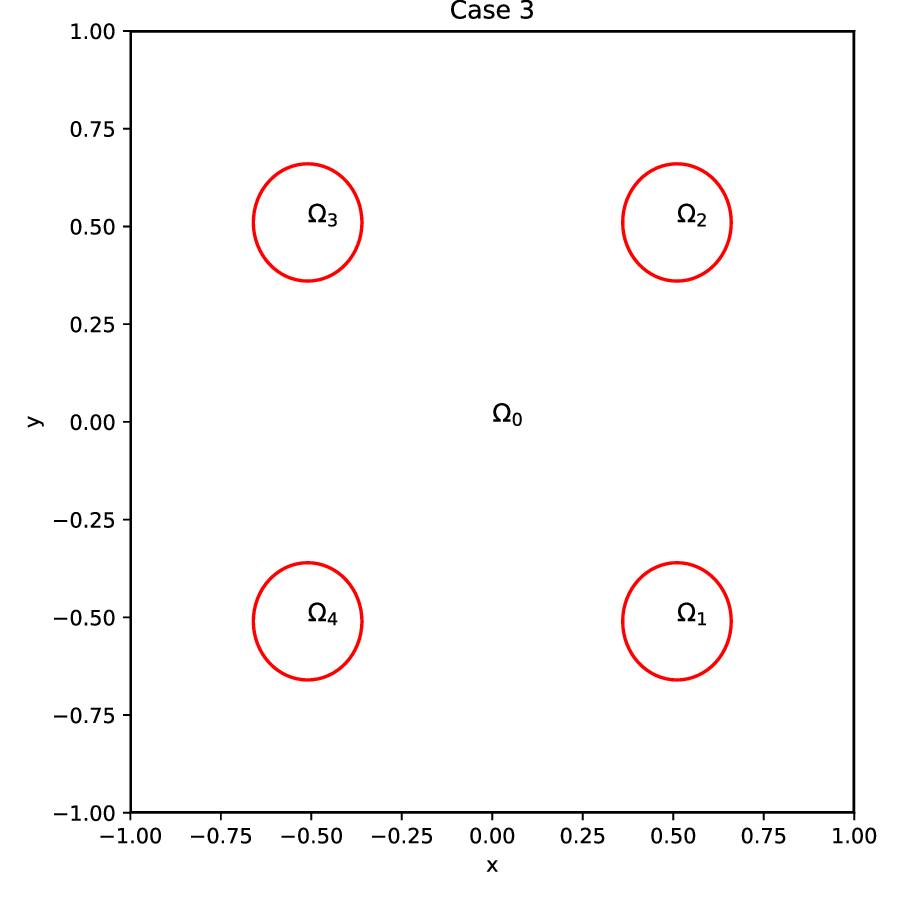}
  \end{minipage}
  \caption{The schematic diagrams of the domain for example \ref{sec:nonlinear-peaks}.}
  \label{Fig:Domain_non}
\end{figure}

\begin{figure}[H]
  \begin{minipage}{0.33\textwidth}
    \centering
    \includegraphics[width=\textwidth, height=5cm]{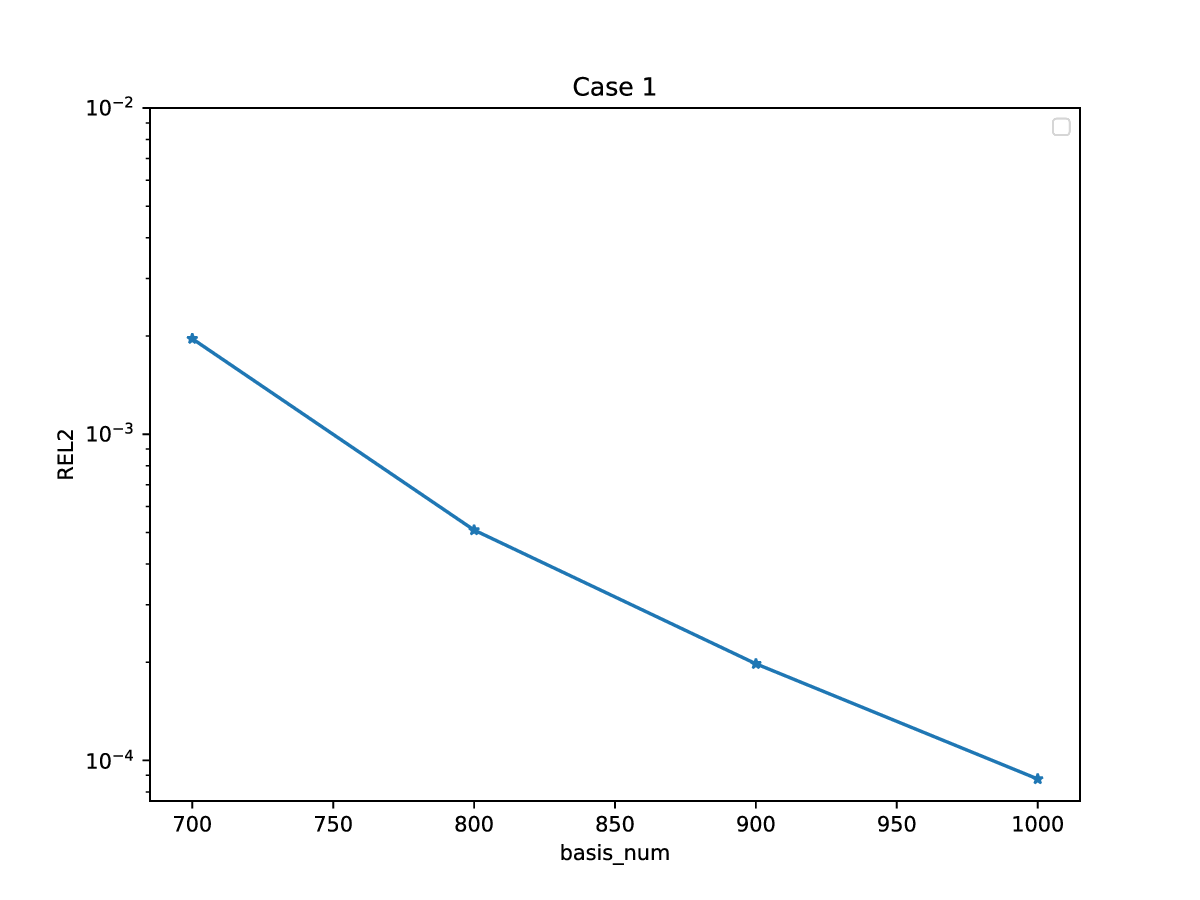}
  \end{minipage}
  \begin{minipage}{0.33\textwidth}
    \centering
    \includegraphics[width=\textwidth, height=5cm]{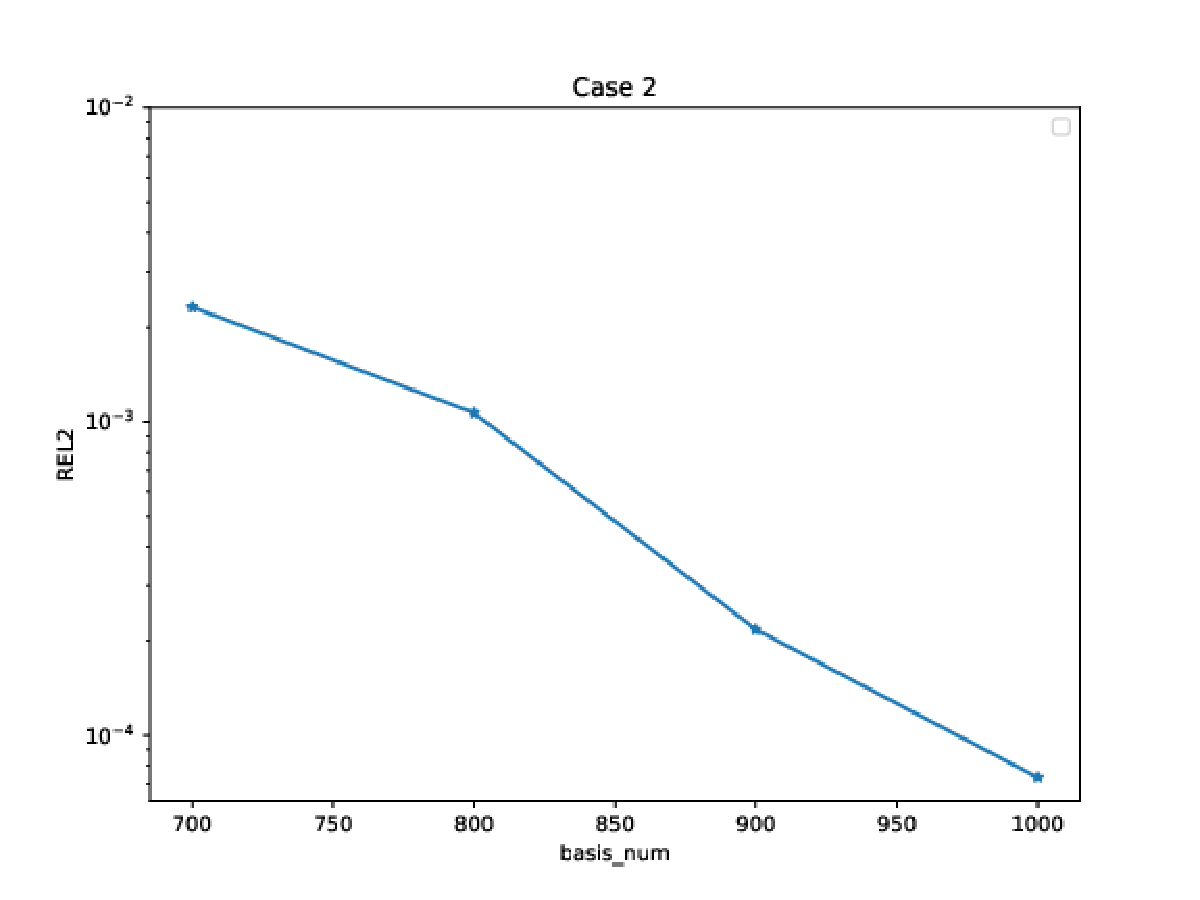}
  \end{minipage}
  \begin{minipage}{0.33\textwidth}
    \centering
    \includegraphics[width=\textwidth, height=5cm]{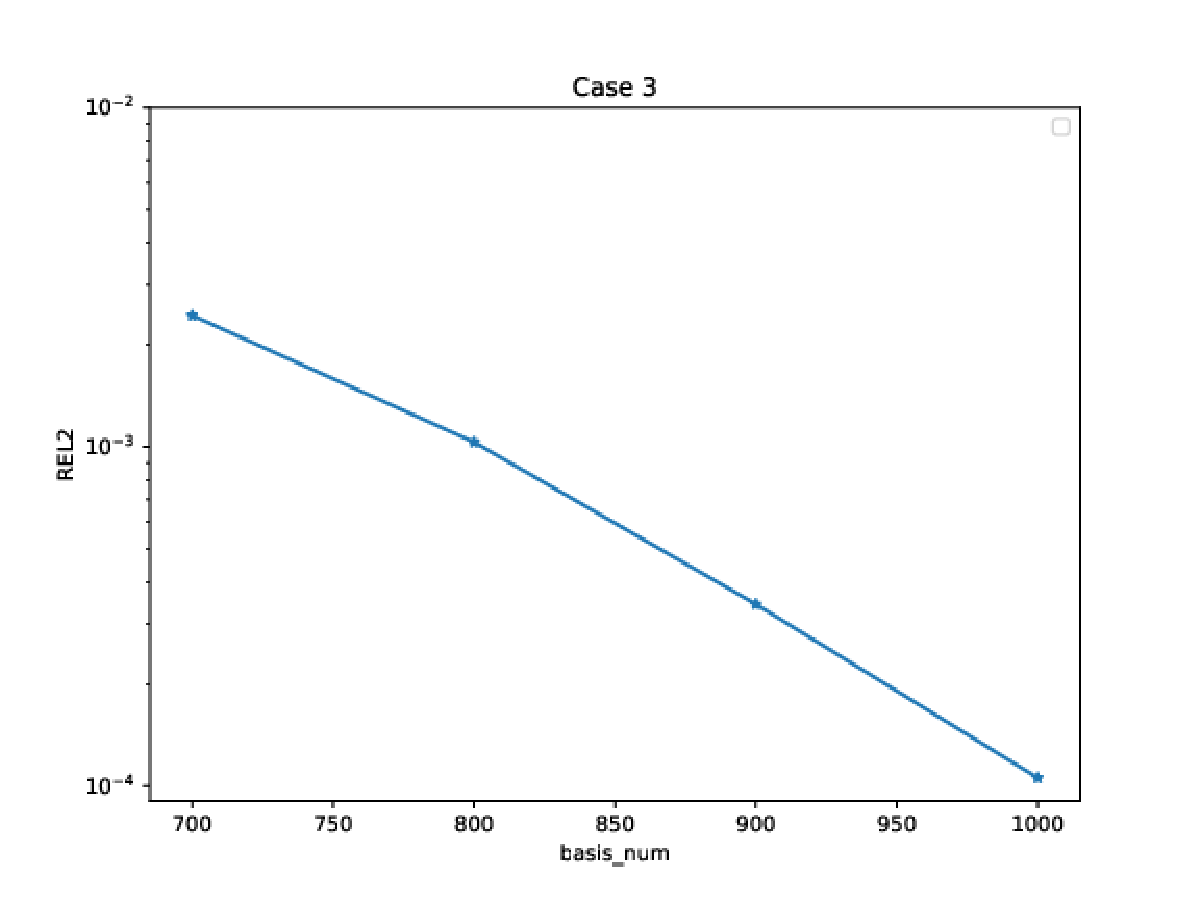}
  \end{minipage}\\
  \caption{The $\operatorname{err}_{L_{2}}$ with the number of scaled neural network basis functions.}
  \label{Fig:rel2_non}
\end{figure}

\begin{figure}[H]
  \begin{minipage}{\textwidth}
    \centering
    \includegraphics[width=\linewidth,height=5cm]{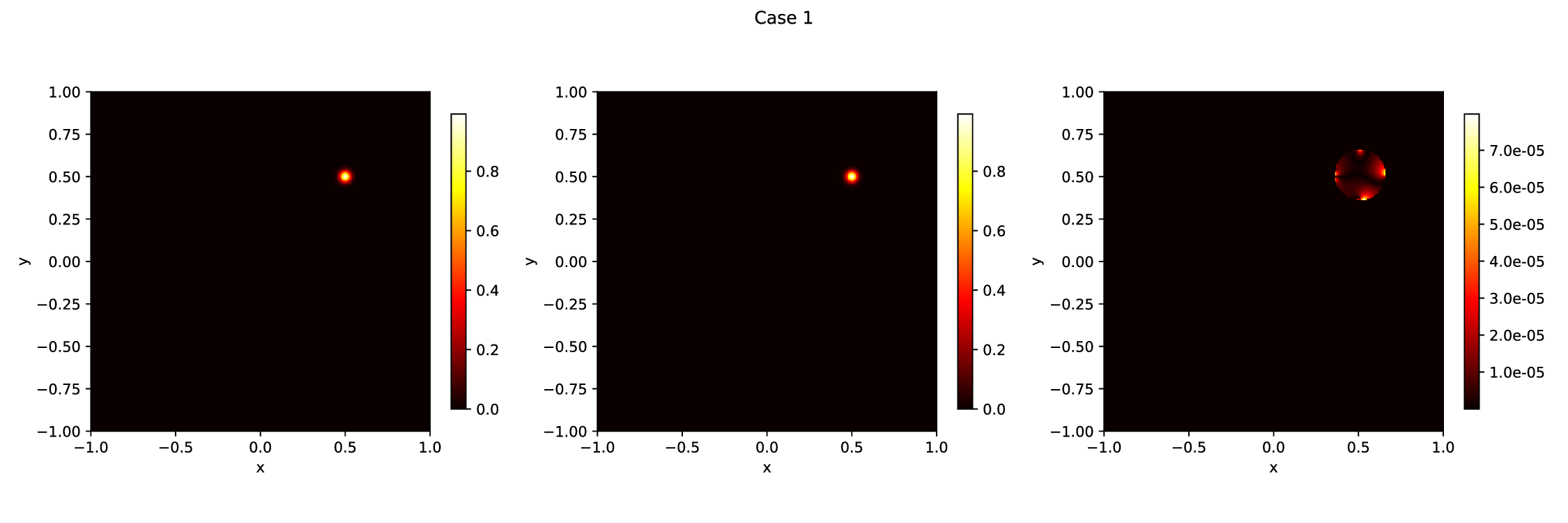}
  \end{minipage}\\
  \begin{minipage}{\textwidth}
    \centering
    \includegraphics[width=\linewidth,height=5cm]{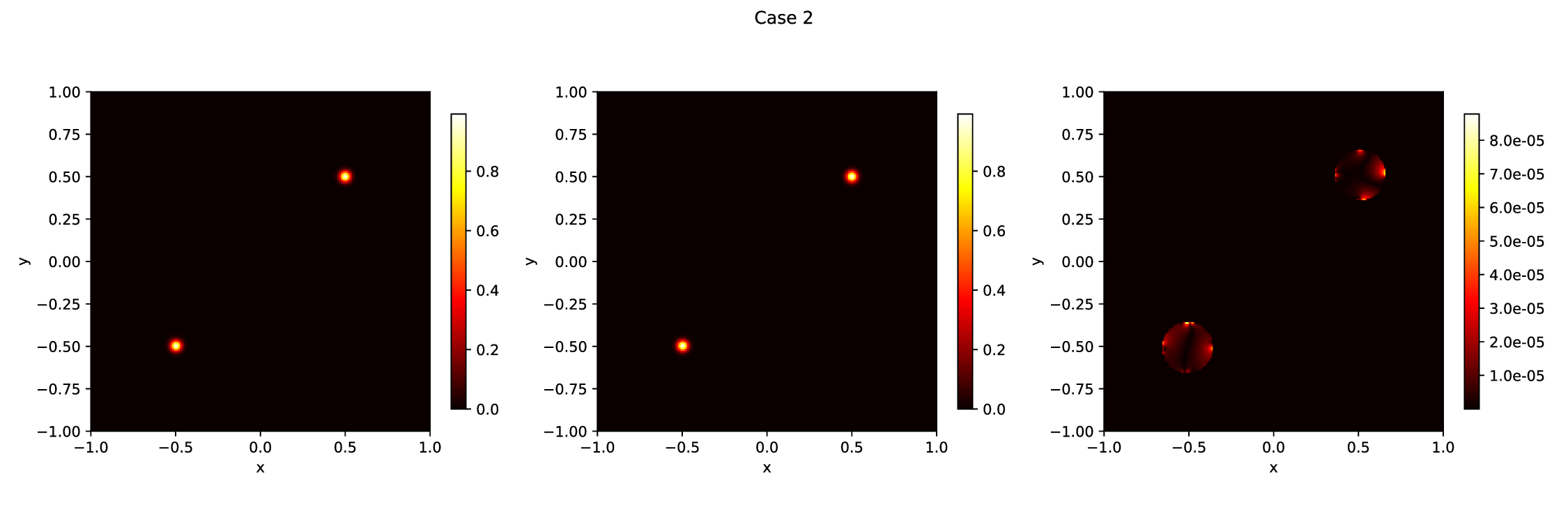}
  \end{minipage}\\
  \begin{minipage}{\textwidth}
    \centering
    \includegraphics[width=\linewidth,height=5cm]{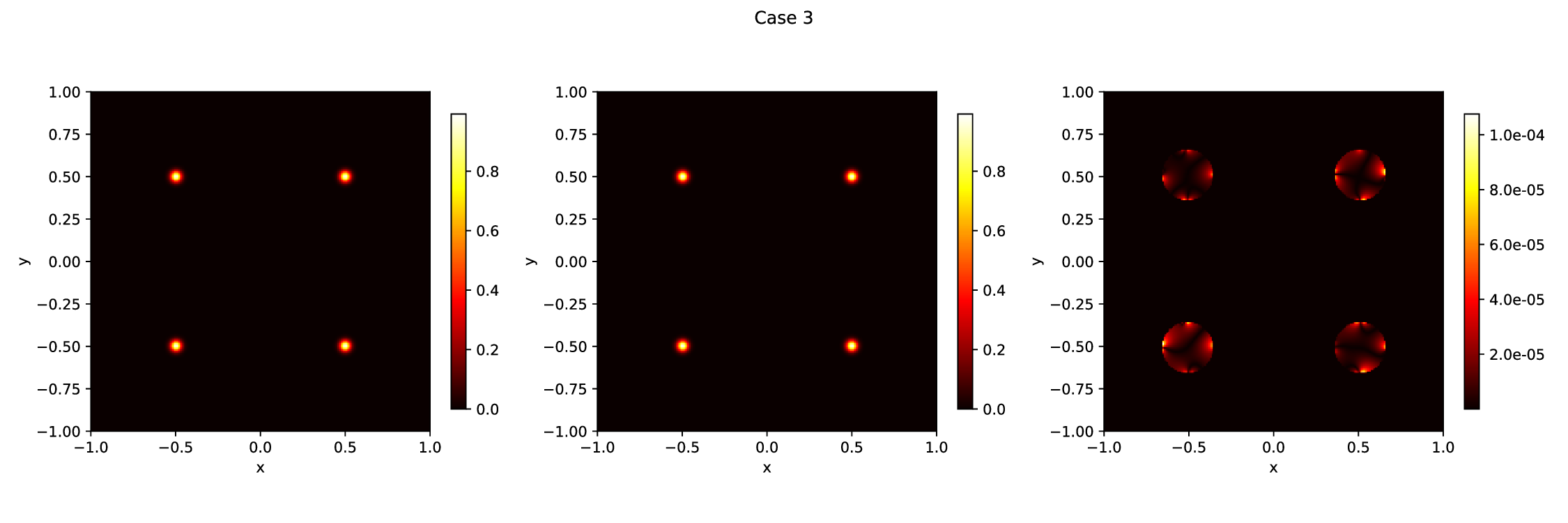}
  \end{minipage}\\
  \caption{The exact solution~(Left), the predicted solution~(Middle) and absolute error~(Right) when $M^{*}=1000$.}
  \label{Fig:heat_cap_non}
\end{figure}

The numerical results are shown in Fig. \ref{Fig:Domain_non},
 \ref{Fig:rel2_non}, \ref{Fig:heat_cap_non} and Tab. \ref{tab:scale-facotr-non}.

From Tab. \ref{tab:scale-facotr-non} and Fig. \ref{Fig:Domain_non}, we can clearly see
that the partition of the domain $\Omega=[-1,1]^2$ and the values of scaling coefficients are the same as the example \ref{sec:linear-2d-peak} for three cases.
 \begin{table}[H]
    \centering
    \begin{tabular}{c|c|c|c}
        & Case 1 & Case 2 & Case 3  \\\hline
        $c_{1}$& 5        & 5        &5 \\\hline
        $c_{2}$& $\times$ & 5        &5 \\\hline
        $c_{3}$& $\times$ & $\times$ &5 \\\hline
        $c_{4}$& $\times$ &$\times$  &5 \\\hline
    \end{tabular}
    \caption{Scale coefficients obtained by Algorithm $\ref{alg:scale-basis}$ for three cases. }
    \label{tab:scale-facotr-non}
\end{table}

Fig. \ref{Fig:rel2_non} shows that the trend of $\operatorname{err}_{L_{2}}$ defined in \eqref{eq:re_l2err}  with the increasing number of the  neural network basis functions $M^{*}$ is the same as the example \ref{sec:linear-2d-peak} for three cases. Besides, we observe that,  the $\operatorname{err}_{L_{2}}$ is smaller than the example \ref{sec:linear-2d-peak} for three cases.
Fig. \ref{Fig:heat_cap_non} presents the predicted solutions and their corresponding absolute error  when $M^{*}=1000$ for three cases. And the performance of absolute error is similar to the example \ref{sec:linear-2d-peak}.


\subsection{Two dimensional problem with corner singularity }\label{sec:corner singularity}
Consider the following two-dimensional Poisson equation:
$$
\begin{aligned}
-\Delta u(x, y)=f(x, y) & & \text { in } \Omega,\\
u(x, y)=g(x, y) & & \text { on } \partial \Omega,
\end{aligned}
$$
where $\Omega$ is $[-1,1]^2\setminus [0,1]^2$ and we specify the true solution as
$$
u(x, y)=(x^2+y^2)^\frac{1}{3}.
$$

In algorithm \ref{alg:ANNB}, let $\epsilon =1e-03,~M_{0}=600,$ and $ r_{k}=r=0.32~(k\ge 1).$  The data set $\bm{X}_{f_{0}}$ consists of a $50 \times 50$ uniform grid points on $[-1,1]^2$, with grid points inside the domain $[0,1]^2$ masked off, i.e., $J_{f_{0}}=1875$. The data set $\bm{X}_{g_{0}}$ is 400 uniformly distributed points on $\partial \Omega,$ i.e., $J_{g_{0}}= 400$. Besides, we compute the $\operatorname{err}_{L_{2}}$ defined in \eqref{eq:re_l2err}  on a $256 \times 256$ uniform grid points  in $[-1,1]^2$ and mask off the grid points outside the domain $[-1,1]^2\setminus[0,1]^2$.

\begin{figure}[H]
    \centering
 \includegraphics[scale=0.35]{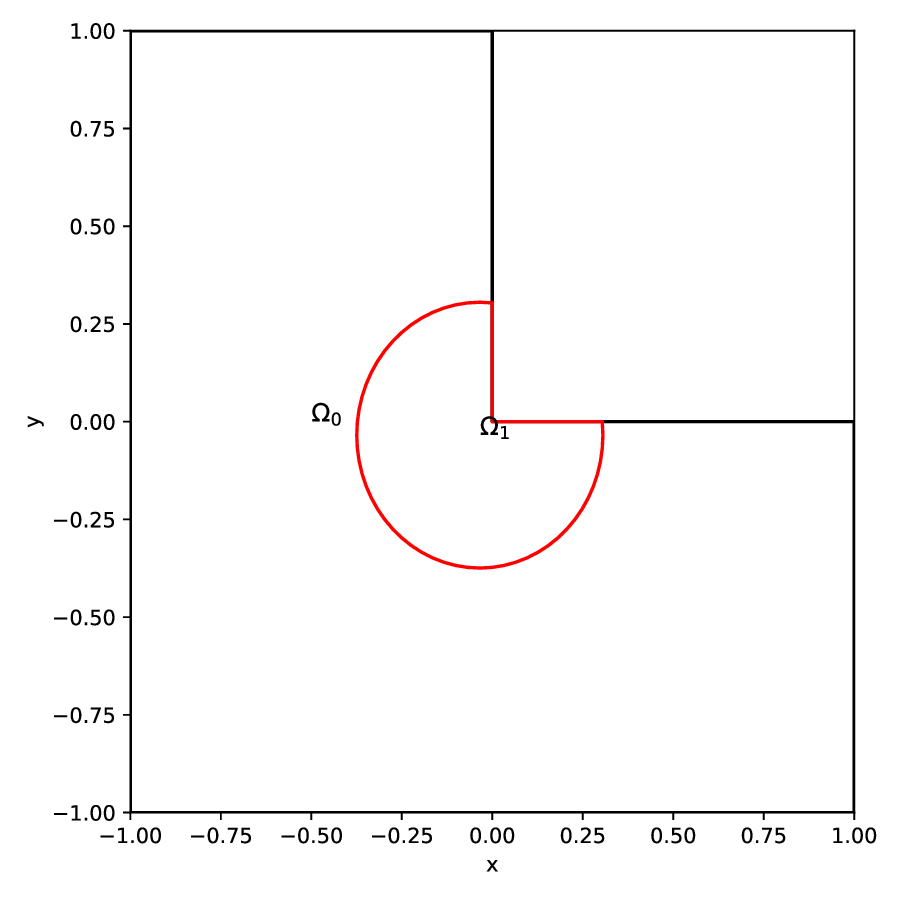}
  \caption{The schematic diagrams of the domain.}
  \label{Fig:Domain_L}
 \end{figure}

\begin{figure}[H]
    \centering
\includegraphics[scale=0.5]{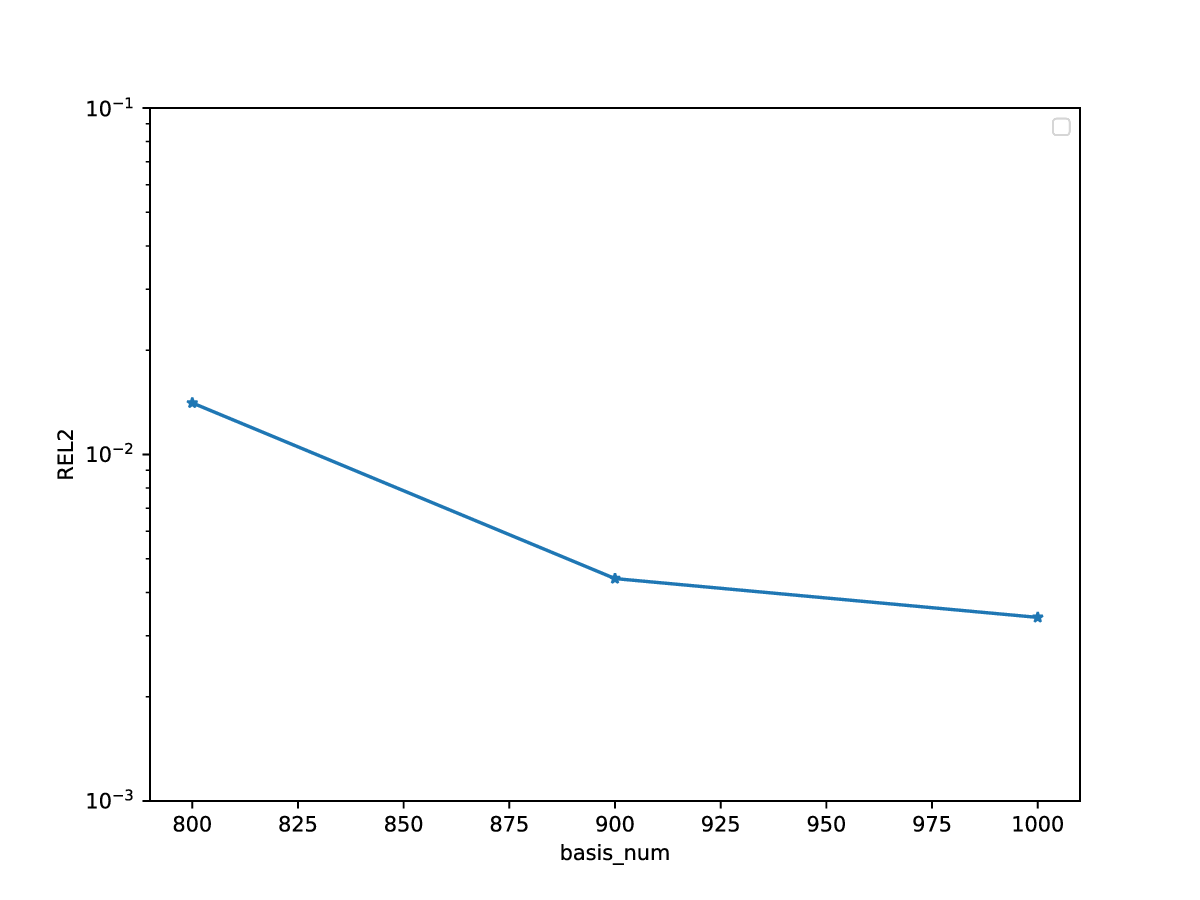}
  \caption{The $\operatorname{err_{L_{2}}}$ with the number of scaled neural network basis functions.}
  \label{Fig:rel2_L}
\end{figure}

\begin{figure}[H]
    \centering
    \includegraphics[width=\textwidth, scale=0.5]{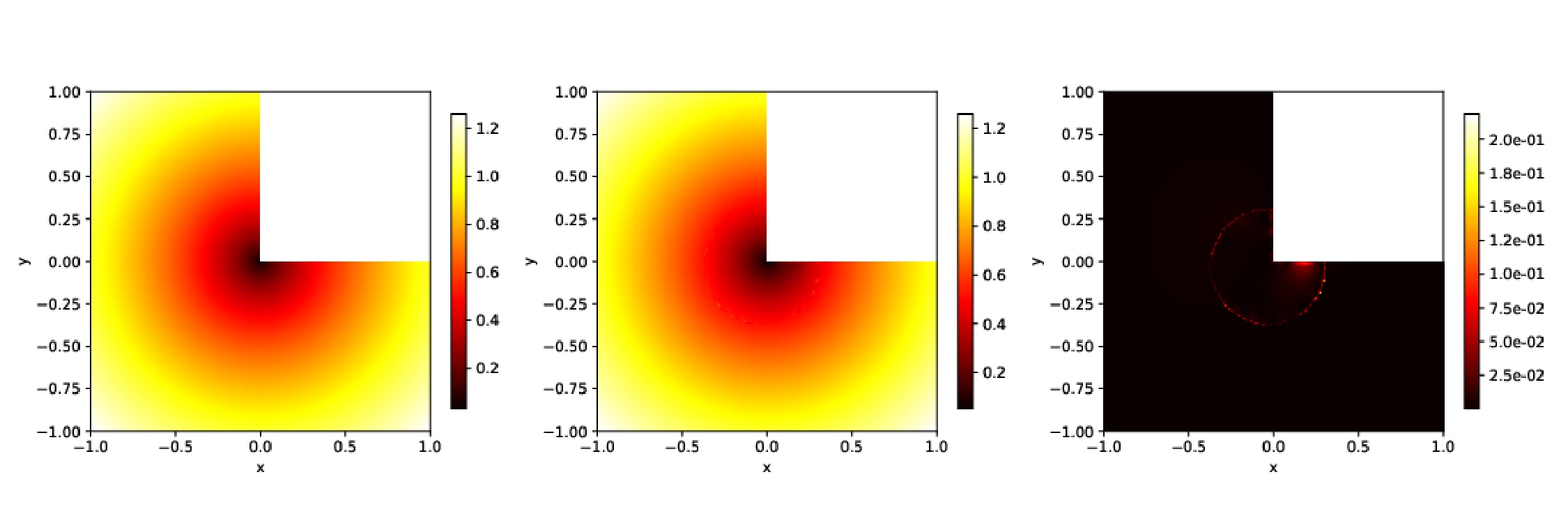}
  \caption{The exact solution~(Left), the predicted solution ~(Middle) and absolute error~(Right) when $M^{*}=1000$.}
  \label{Fig:heat_cap_L}
\end{figure}

The numerical results are presented in Fig. \ref{Fig:Domain_L}, \ref{Fig:rel2_L}, and \ref{Fig:heat_cap_L}. Fig. \ref{Fig:Domain_L} illustrates the partitioning of the domain \( \Omega \) into two sub-domains \( \Omega_{0} \) and \( \Omega_{1} \), where \( \Omega_{1} = B_{\bm{x}_{1}}(r) \) with \( \bm{x}_{1} = (-0.0345, -0.0345)^{\transpose} \). Additionally, the scaling coefficient \( c_{1} \), obtained using algorithm \ref{alg:scale-basis}, equals 5 when \( M^{*} = 800, 900 \) and equals 4 when \( M^{*} = 1000 \).

Fig. \ref{Fig:rel2_L} shows that  $\operatorname{err}_{L_{2}}$  decreases as the number of neural network basis functions $M^{*}$ increases. Specifically,  $M^{*}$ is chosen as $800,900,1000$ and the $\operatorname{err}_{L_{2}}$ reaches approximately  $1\mathrm{e}-3$ when $M^*=1000$. Fig. \ref{Fig:heat_cap_L} presents the predicted solutions and their corresponding absolute error when $M^{*}=1000$.

\subsection{Three-dimensional problem.}\label{sec:ex-3d}
Consider the following three dimensional Poisson equation:
$$
\begin{aligned}
-\Delta u(x, y,z)=f(x, y,z) & & \text { in } \Omega, \\
u(x, y,z)=g(x, y,z) & & \text { on } \partial \Omega,
\end{aligned}
$$
where $\Omega$ is $[-1,1]^3$ and we specify the true solution as
$$
u(x, y)=\exp \left(-1000\left[(x-0.5)^2+(y-0.5)^2+(z-0.5)^2 \right]\right).
$$
 Besides, we compute the  $\operatorname{err}_{L_{2}}$ defined in  \eqref{eq:re_l2err} on a $50 \times 50\times50 $ uniform mesh in $\Omega=[-1,1]^3$.

 In Algorithm \ref{alg:ANNB}, let $\epsilon =1e-04,~M_{0}=2000$ and $r_{k} = r = 0.11~(k\ge 1).$  The data set $\bm{X}_{f_{0}}$ consists of 10000 uniformly distributed points in $\Omega=[-1, 1]^3$, i.e., $J_{f_{0}}=10000$ and the data set $\bm{X}_{g_{0}}$ consists of  $2400$ uniformly distributed points on $\partial\Omega$, 400 points on each side of $\partial\Omega$, i.e., $J_{g_{0}}=2400$.

\begin{figure}[H]
    \centering    \includegraphics[scale=0.3]{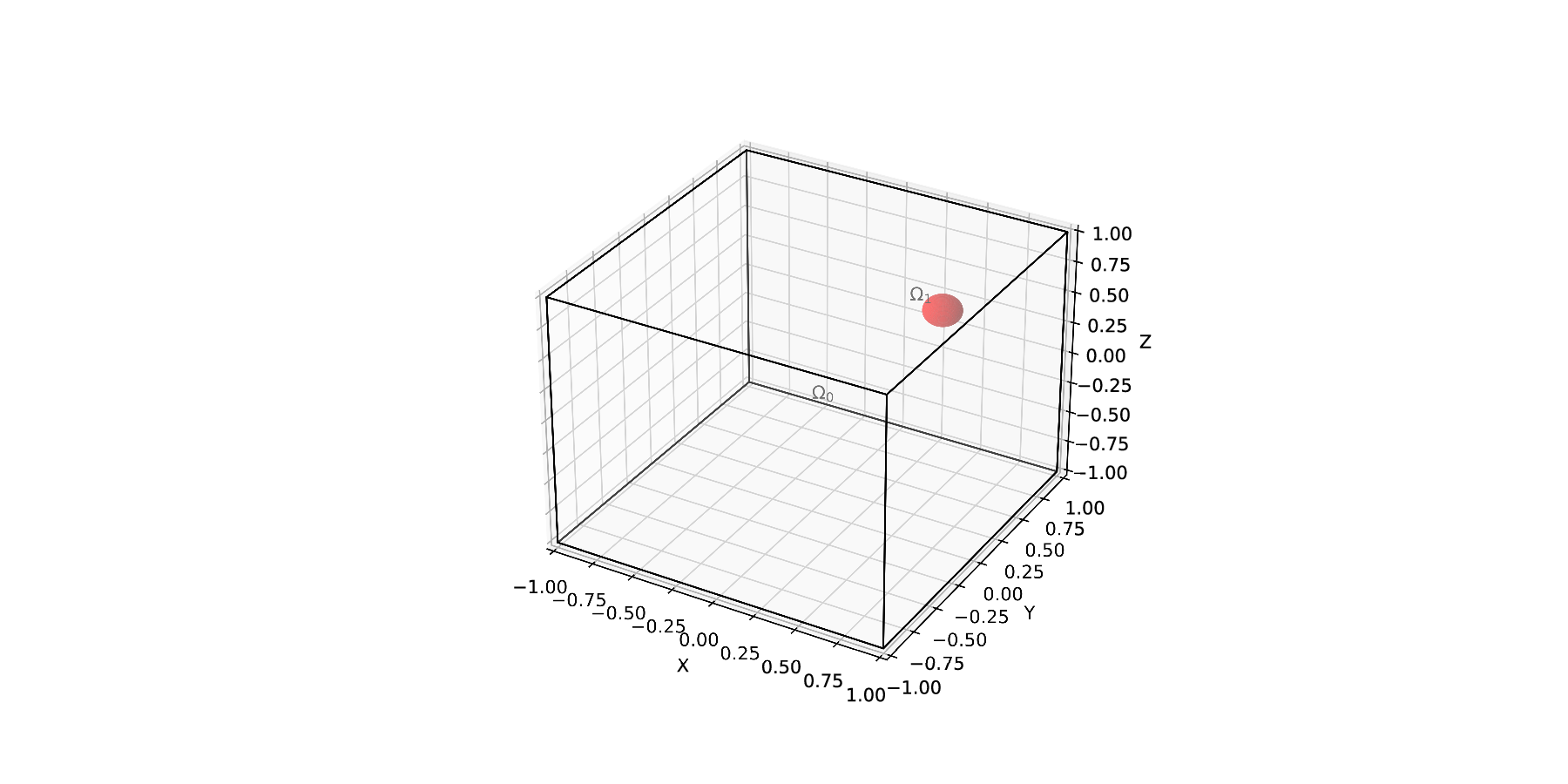}
    \caption{\small The schematic diagrams of the domain.}
    \label{Fig:domain_3d}
\end{figure}

\begin{figure}[H]
    \centering
\includegraphics[scale=0.5]{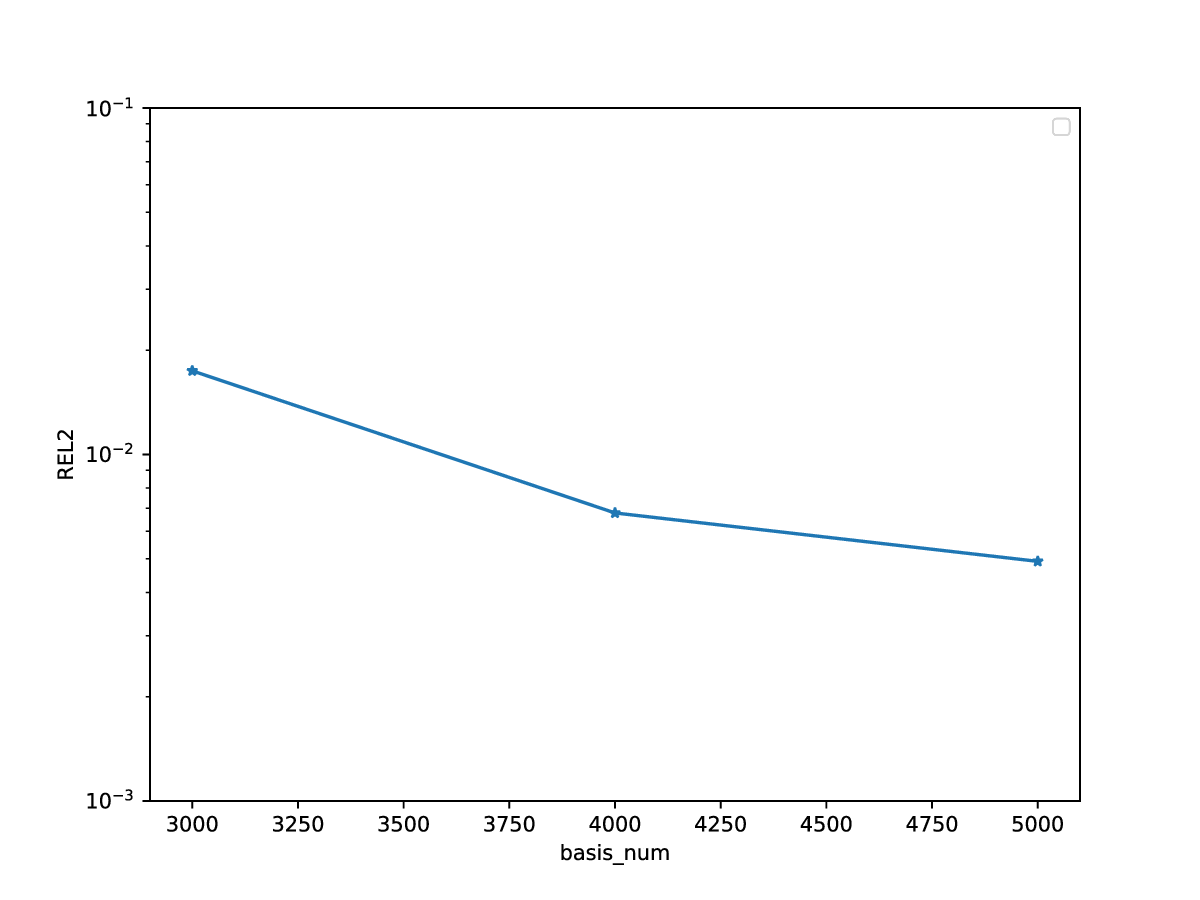}
    \caption{\small The $\operatorname{err}_{L_{2}}$  with the number of scaled neural network bases.}
    \label{Fig:rel2_3d}
    \end{figure}

\begin{figure}[H]
    \centering
\includegraphics[scale=0.4]{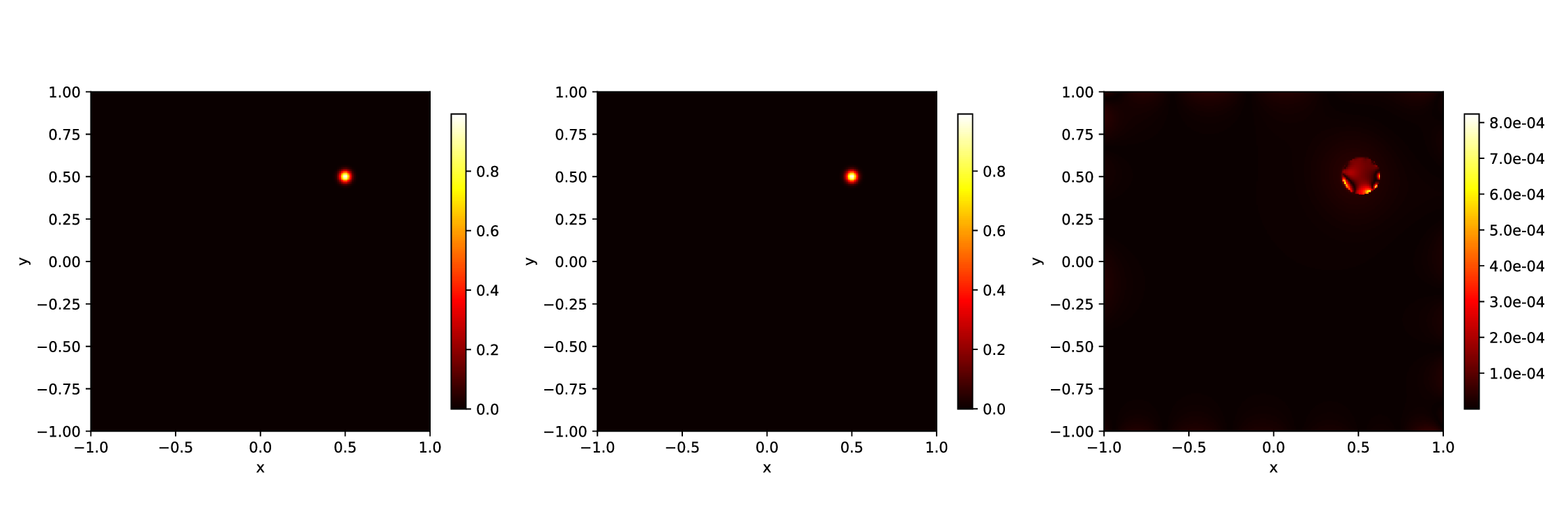}
\caption{\small The exact solution~(Left), the predicted solution~(Middle) and absolute error~(Right) at $z=0.5$ when $M^*=5000$ .}
  \label{Fig:heat_cap_3d}
\end{figure}

The numerical results are shown in Fig. \ref{Fig:domain_3d}, \ref{Fig:rel2_3d} and \ref{Fig:heat_cap_3d}. Fig. \ref{Fig:domain_3d} shows that the domain $\Omega$ is eventually partitioned into two domains $\Omega_{0}$ and $\Omega_{1}$, where $\Omega_{1}=B_{\bm{x}_{1}}(r)$ with $\bm{x}_{1}=(0.5170,0.5050,0.5000)^{\transpose}$.
Besides, the scaling coefficient $c_{1}$ obtained by algorithm \ref{alg:scale-basis} equals 6 when $M^{*}=3000$, equals 7 when $M^{*}=4000$ and equals 8 when $M^{*}=5000$.

It is clearly seen from Fig. \ref{Fig:rel2_3d}  that the $\operatorname{err}_{L_{2}}$ defined in \eqref{eq:re_l2err}  decreases  as  the number of neural network basis functions  $M^{*}$ increases and $\operatorname{err}_{L_{2}}$  reaches approximately $1\mathrm{e}-3$ when $M^{*}=5000$.

Fig. \ref{Fig:heat_cap_3d} presents the predicted solutions and their corresponding absolute errors at $z=0.5$ when $M^{*}=5000$.
\section{Conclusion}
In this paper, we have proposed the ANNB method for partial differential equations with low-regular solutions. Then, we apply the method to numerically solve some peak problems; the numerical performance is efficient. The work can be viewed as our first attempt to attack such class of problems using the adaptive method combined with the neural network basis method. There are quite a number of issues worth further study. For instance, how to find a way to figure out all peak domains in a single iteration, how to find a better method in determining the scaling coefficients, and how to handle time-dependent problems with low-regular solutions.
\section*{Appendix}
In this appendix, we present the procedure for solving problem \eqref{pro:min-scale factor-s}. When the differential operator \( \mathcal{L} \) is linear, the problem \eqref{pro:min-scale factor-s} can be rewritten as the following linear least squares problem:
\begin{equation}\label{pro:K-min-linear}
  \begin{aligned}
  \min_{\bm{\alpha}^{0}, \cdots, \bm{\alpha}^{K}} &\left\{
\left(\sum_{\bm{x}_{f_{K}}\in \bm{X}_{f_{K}}}\left|\sum_{m=0}^{M_{K}} \alpha_{m,K}\mathcal{L}\psi_{m,K}(\bm{x}_{f_{K}},s)-f(\bm{x}_{f_{K}})\right|^2\right.\right.\\
  &\left. +\sum_{\bm{x}_{g_{K}}\in \bm{X}_{g_{K}}}\left|\sum_{m=0}^{M_{K}} \alpha_{m,K}\psi_{m,K}(\bm{x}_{g_{K}},s)-g(\bm{x}_{g_{K}})\right|^2\right)\\
  &+\sum_{\bm{x}_{\Gamma_{K}}\in \bm{X}_{\Gamma_{K}}}\left(\left| \sum_{m=0}^{M_{K}} \alpha_{m,K} \psi_{m,K}(\bm{x}_{\Gamma_{K}},s)- \sum_{m=0}^{M_{0}} \alpha_{m,0} \psi_{m,K}(\bm{x}_{\Gamma_{K}},s)\right|^2\right.\\
  &\left.\left.+\left|\frac{\sum_{m=0}^{M_{K}}\alpha_{m,K}\partial \psi_{m,K}(\bm{x}_{\Gamma_{K}},s)}{\partial \bm{n}_{K}}-\frac{\sum_{m=0}^{M_{0}}\partial\alpha_{m,0}\ \psi_{m,0}(\bm{x}_{\Gamma_{K}},s)}{\partial \bm{n}_{K}}\right|^2\right)\right\}
  \end{aligned}
  \end{equation}

When the differential operator \( \mathcal{L} \) is nonlinear, the Gauss-Newton method is used to solve the least squares problem \eqref{pro:min-scale factor-s}. At the \( n \)-th iteration step, we denote \( u_{K}^{n}(\bm{x}) \) as the approximation of the solution \( u_{K}(\bm{x}) \), and \( v_{K}^{n}(\bm{x}) \) as the increment field that needs to be computed. Then,
$$
u_{K}^{n+1}(\bm{x})= u_{K}^{n}(\bm{x})+v_{K}^{n}(\bm{x}).
$$


At the Newton step \( n \), we represent \( u_{K}^{n}(\bm{x}) \) and \( v_{K}^{n}(\bm{x}) \) as follows.
$$
u_{K}^{n}(x)\approx\Tilde{u}_{K}(\bm{\alpha}_{K}^{n},\bm{x}) =\sum_{m=0}^{M_{K}} \alpha_{m,K}^{n} \psi_{m,K}(\bm{x},s),~
v_{K}^{n}(x)\approx\Tilde{u}_{K}(\bm{a}_{K}^{n},\bm{x})= \sum_{m=0}^{M_{K}} a_{m,K}^{n} \psi_{m,K}(\bm{x},s),
$$
where \( \bm{\alpha}_{K}^{n}=(\alpha_{0,K}^{n},\alpha_{1,K}^{n},\cdots,\alpha_{M_{K},n}^{K})^{\top} \) and \( \bm{a}_{K}^{n}=(a^{K}_{0,n},a^{K}_{1,n},\cdots,a^{K}_{M_{K},n})^\top \). Then, the unknown parameters \( \bm{a}_{K} \) can be obtained by solving the following linear least squares problem:
\begin{equation}\label{pro:K-Newton-min}
  \begin{aligned}
    \min_{\bm{a}_{K}^{n}} &\left\{
    \left(\sum_{\bm{x}_{f_{K}}\in \bm{X}_{f_{K}}}\left|\sum_{m=0}^{M_{K}} a_{m,K}^{n}D\mathcal{L}(\tilde{u}_{0}^{n};\psi_{m,K})(\bm{x}_{f_{K}})-f(\bm{x}_{f_{K}})+\mathcal{L}\tilde{u}_{0}^{n}(\bm{x}_{f_{K}})\right|^2 \right.\right.\\
    \qquad \qquad &\left.+\sum_{\bm{x}_{g_{K}}\in \bm{X}_{g_{K}}}\left|\sum_{m=0}^{M_{K}} a_{m,K}^{n}\psi_{m,K}(\bm{x}_{g_{K}})-g(\bm{x}_{g_{K}})+\tilde{u}_{K}^{n}(\bm{x}_{g_{K}})\right|^2\right)\\
    &+\sum_{\bm{x}_{\Gamma_{K}}\in \bm{X}_{\Gamma_{K}}}\left(\left| \sum_{m=0}^{M_{K}} a_{m,K} \psi_{m,K}(\bm{x}_{\Gamma_{K}})- \sum_{m=0}^{M_{0}} a_{m,0} \psi_{m,K}(\bm{x}_{\Gamma_{K}})+\tilde{u}_{K}(\bm{\alpha}_{K}^{n},\bm{x}_{\Gamma_{K}})-\tilde{u}_{0}(\bm{\alpha}_{0}^{n},\bm{x}_{\Gamma_{K}})\right|^2\right.\\
    &\left.\left.+\left|\frac{\sum_{m=0}^{M_{K}}\alpha_{m,K}\partial \psi_{m,K}(\bm{x}_{\Gamma_{K}})}{\partial \bm{n}_{K}}-\frac{\sum_{m=0}^{M_{0}}\alpha_{m,0}\ \psi_{m,0}(\bm{x}_{\Gamma_{K}})}{\partial \bm{n}_{K}}+\frac{\partial\tilde{u}_{K}(\bm{\alpha}_{K}^{n},\bm{x}_{\Gamma_{K}})}{\partial\bm{n}_{K}}-\frac{\partial\tilde{u}_{0}(\bm{\alpha}_{0}^{n},\bm{x}_{\Gamma_{K}})}{\partial\bm{n}_{K}}\right|^2\right)\right\},
    \end{aligned}
\end{equation}

which can be solved directly.

\bibliographystyle{abbrv}
\bibliography{ref1}
\end{document}